\documentclass[11pt,a4paper]{article}

\usepackage{theorem,enumerate}
\usepackage{amsmath,latexsym,amssymb,amsfonts}
\usepackage{graphicx}

\newtheorem{thm}{Theorem}
\newtheorem{Thm}{Theorem}

\newtheorem{cor}[thm]{Corollary}
\newtheorem{lem}{Lemma}
\newtheorem{prop}[thm]{Proposition}

\newtheorem{rem}{Remark}

\newtheorem{Conj}[Thm]{Conjecture}
\newtheorem{claim}{Claim}
\newtheorem{subclaim}{Subclaim}[claim]

\newtheorem{problem}{Problem}

\newcommand{\proof}{\medbreak\noindent\textit{Proof.}\quad}

\newcommand{\qed}{{$\quad\square$\vs{3.6}}}

\newcommand{\vs}[1]{\vspace*{#1 mm}}

\numberwithin{equation}{section}

\newcounter{mymemory2}

\title{On the existence of vertex-disjoint subgraphs\\ with high degree sum\thanks{An extended abstract has been accepted in: 
EuroComb 2015, Electr. Notes Discrete Math., vol. 49, 2015, pp. 359--366.}}

\author{Shuya Chiba$^{1}$\thanks{This work was supported by JSPS KAKENHI Grant Number 26800083}
\thanks{E-mail address: \texttt{schiba@kumamoto-u.ac.jp}} 
\hspace{+12pt}
Nicolas Lichiardopol$^{2}$\thanks{E-mail address: \texttt{nicolas.lichiardopol@neuf.fr}}
\vspace{+8pt}
 \\
\small
$^1$\small\textsl{Applied Mathematics, 
Faculty of Advanced Science and Technology, 
Kumamoto University,}\\ 
\vspace{+6pt}
\small\textsl{2-39-1 Kurokami, Kumamoto 860-8555, Japan}\\
\small
$^{2}$\small\textsl{Lyc\'{e}e A. de Craponne, Salon, France}}

\date{}

\begin{document}
\setlength{\baselineskip}{17pt}

\maketitle

\vspace{-24pt}
\begin{abstract} 
For a graph $G$, 
we denote by $\sigma_{2}(G)$ the minimum degree sum of two non-adjacent vertices 
if $G$ is non-complete; otherwise, $\sigma_{2}(G) = +\infty$. 
In this paper, 
we prove the following two results: 
(i) If $s_{1}, s_{2} \ge 2$ are integers 
and $G$ is a non-complete graph with $\sigma_{2}(G) \ge 2(s_{1} + s_{2} + 1) - 1$, 
then $G$ contains two vertex-disjoint subgraphs $H_{1}$ and $H_{2}$ such that 
each $H_{i}$ is a graph of order at least $s_{i}+1$ with $\sigma_{2}(H_{i}) \ge 2s_{i} - 1$. 
(ii) If $s_{1}, s_{2} \ge 2$ are integers 
and $G$ is a triangle-free graph of order at least $3$ with $\sigma_{2}(G) \ge 2(s_{1} + s_{2}) - 1$, 
then $G$ contains two vertex-disjoint subgraphs $H_{1}$ and $H_{2}$ such that 
each $H_{i}$ is a graph of order at least $2s_{i}$ with $\sigma_{2}(H_{i}) \ge 2s_{i} - 1$. 
By using this result, we also give some corollaries 
concerning degree conditions for the existence of $k$ vertex-disjoint cycles.



\end{abstract}

\section{Introduction}
\label{introduction} 

In this paper, we consider finite simple graphs, 
which have neither loops nor multiple edges.
For terminology and notation not defined in this paper, we refer the readers to \cite{Diestel}. 
Let $G$ be a graph. 
We denote by $V(G)$, $E(G)$ and $\delta(G)$ 
the vertex set, the edge set and the minimum degree of $G$, respectively. 
We write $|G|$ for the order of $G$, that is, $|G|=|V(G)|$. 
We denote by $d_{G}(v)$ the degree of a vertex $v$ in $G$. 
If $H$ is a subgraph of $G$, 
then $d_{H}(v)$ is the number of vertices in $H$ that are adjacent to a vertex $v$ of $G$. 
The invariant $\sigma_2(G)$ is defined to be the minimum degree sum
of two non-adjacent vertices of $G$, i.e., 
$\sigma_2(G)=\min \big\{ d_G(u) + d_G(v): u, v \in V(G), u \neq v, uv \notin E(G) \big\}$
if $G$ is non-complete; 
otherwise, 
let $\sigma_2(G) = + \infty$. 
We denote by $g(G)$ the \textit{girth} of $G$, 
i.e., the length of a shortest cycle of $G$. 
In this paper, ``disjoint'' always means ``vertex-disjoint''. 
A pair $(H_{1}, H_{2})$ is called a \textit{partition of} $G$ 
if $H_{1}$ and $H_{2}$ are two disjoint induced subgraphs of $G$ 
such that $V(G) = V(H_{1}) \cup V(H_{2})$.

\medskip

Stiebitz \cite{Stiebitz} considered the decomposition of graphs under degree constraints 
and proved the following result.

\begin{Thm}[Stiebitz \cite{Stiebitz}]
\label{Stiebitz}
Let $s_{1}, s_{2} \ge 1$ be integers, 
and let $G$ be a graph. 
If $\delta(G) \ge s_{1} + s_{2} + 1$, 
then there exists a partition $(H_{1}, H_{2})$ of $G$ such that 
$\delta(H_{i}) \ge s_{i}$ for $i \in \{ 1, 2 \}$. 
\end{Thm}

Kaneko \cite{Kaneko} showed that 
the same holds for triangle-free graphs with minimum degree at least $s_{1} + s_{2}$.

\begin{Thm}[Kaneko \cite{Kaneko}]
\label{Kaneko}
Let $s_{1}, s_{2} \ge 1$ be integers, 
and let $G$ be a graph. 
If $\delta(G) \ge s_{1} + s_{2}$ and $g(G) \ge 4$, 
then there exists a partition $(H_{1}, H_{2})$ of $G$ such that 
$\delta(H_{i}) \ge s_{i}$ for $i \in \{ 1, 2 \}$. 
\end{Thm}

Diwan further improved Theorem~\ref{Stiebitz} for graphs with girth at least $5$, see \cite{Diwan}. 
Bazgan, Tuza and Vanderpooten \cite{Bazgan et al.} gave polynomial-time algorithms that find such partitions.

The purpose of this paper is to consider $\sigma_{2}$-versions of Theorems~\ref{Stiebitz} and \ref{Kaneko}. 
More precisely, 
we consider the following problems.

\begin{problem} 
\label{partition problem sigma2}
Let $s_{1}, s_{2} \ge 2$ be integers, 
and let $G$ be a non-complete graph. 
If $\sigma_{2}(G) \ge 2(s_{1} + s_{2} + 1) - 1$, 
determine whether 
there exists a partition $(H_{1}, H_{2})$ of $G$ such that 
$\sigma_{2}(H_{i}) \ge 2s_{i} - 1$ 
and $|H_{i}| \ge s_{i} + 1$ for $i \in \{ 1, 2 \}$. 
\end{problem}

\begin{problem} 
\label{partition problem sigma2 tri-free}
Let $s_{1}, s_{2} \ge 2$ be integers, 
and let $G$ be a graph of order at least $3$. 
If $\sigma_{2}(G) \ge 2(s_{1} + s_{2}) - 1$ and $g(G) \ge 4$, 
determine whether 
there exists a partition $(H_{1}, H_{2})$ of $G$ such that 
$\sigma_{2}(H_{i}) \ge 2s_{i} - 1$ 
and 
$|H_{i}| \ge 2s_{i}$ for $i \in \{ 1, 2 \}$. 
\end{problem}

In Problem~\ref{partition problem sigma2} (resp., Problem~\ref{partition problem sigma2 tri-free}), 
if we drop the condition ``$|H_{i}| \ge s_{i} + 1$ (resp., $|H_{i}| \ge 2s_{i}$)'' in the conclusion, 
then it is an easy problem. 
Because, 
for each edge $xy$ in a graph $G$ satisfying the assumption of Problem~\ref{partition problem sigma2} 
(resp., the assumption of Problem~\ref{partition problem sigma2 tri-free}), 
$H_{1} = G[\{x, y\}]$ and $H_{2} = G - \{x, y\}$ satisfy 
$\sigma_{2}(H_{1}) = \infty > 2s_{1} - 1$ and $\sigma_{2}(H_{2}) \ge \sigma_{2}(G) - 2|\{x, y\}| \ge 2s_{2} - 1$. 
Here, for a vertex subset $X$ of a graph $G$, 
$G[X]$ denotes the subgraph of $G$ induced by $X$, 
and let $G - X = G[V(G) \setminus X]$. 
(Similarly, for the case where $s_{i} = 1$ for some $i$, we can easily solve it.)

If $G_{1}$ is a balanced complete multipartite graph with $r+1 \ ( \ge 4)$ partite sets of size $s \ ( \ge 2)$, 
then $\sigma_{2}(G_{1}) = 2rs = 2 \big( (rs- r + 1) + (r-1) + 1 \big) - 2$, and 
we can check that $G_{1}$ contains no partitions as in Problem~\ref{partition problem sigma2} for $(s_{1}, s_{2}) = (rs- r + 1, r -1)$. 
Thus, the condition ``$\sigma_{2}(G) \ge 2(s_{1} + s_{2}+1) - 1$'' in Problem~\ref{partition problem sigma2} is best possible in a sense 
if it's true. 
If $G_{2}$ is a complete bipartite graph $K_{s_{1} + s_{2} - 1, s_{1}+s_{2}}$, 
then $\sigma_{2}(G_{2}) = 2(s_{1} + s_{2}) - 2$, 
and $G_{2}$ does not contain partitions as in Problem~\ref{partition problem sigma2 tri-free}. 
Thus, $G_{2}$ shows that 
the condition ``$\sigma_{2}(G) \ge 2(s_{1} + s_{2}) - 1$'' 
in Problem~\ref{partition problem sigma2 tri-free} 
is also best possible if it's true.

\medskip

Before giving the main result, we introduce the outline of the proof of Theorems~\ref{Stiebitz} and \ref{Kaneko}. 
The proof consists of the following two steps: 
\begin{enumerate}[{\textup{{\bf Step \arabic{enumi}:}}}]
\setlength{\parskip}{1.5pt}
\setlength{\itemsep}{1.5pt}
\item
To show the existence of two disjoint subgraphs of high minimum degree, 
i.e., 
we show the existence of two disjoint subgraphs $H_{1}$ and $H_{2}$ such that $\delta(H_{i}) \ge s_{i}$ for $i \in \{1, 2\}$. 
\item
To show the existence of two disjoint subgraphs of high minimum degree that partition $V(G)$ 
by using Step~1. 
\end{enumerate}
In particular, 
in the proof of Theorems~\ref{Stiebitz} and \ref{Kaneko}, 
Step~2 follows easily from Step~1. 
In fact, 
if $G$ is a graph with $\delta(G) \ge s_{1} + s_{2} - 1$, 
and $G$ contains 
a pair $(H_{1}, H_{2})$ of disjoint subgraphs with $\delta(H_{i}) \ge s_{i}$ for $i \in \{1, 2\}$, 
then we can easily transform the pair into a partition of $G$ keeping its minimum degree condition 
(see \cite[Proposition 4]{Stiebitz}).

Considering the situation for the proof of Theorems~\ref{Stiebitz} and \ref{Kaneko}, 
one may approach Problems~\ref{partition problem sigma2} and~\ref{partition problem sigma2 tri-free} 
by following the same steps as above. 
However, 
for the case of $\sigma_{2}$-versions, 
neither Step~1 nor Step~2 is an easy problem 
because we allow vertices with low degree. 
In fact, in the proof of Step~2 for Theorem~\ref{Stiebitz} (\cite[Proposition~4]{Stiebitz}), 
the assumption that every vertex has high degree plays a crucial role. 
At the moment, 
we don't know 
whether we can extend disjoint subgraphs of high minimum ``degree sum'' to a partition or not. 
However, 
we can solve Step~1 for Problems~\ref{partition problem sigma2} and \ref{partition problem sigma2 tri-free}. 
The following are our main results.

\begin{thm} 
\label{disjoint subgraphs sigma2}
Let $s_{1}, s_{2} \ge 2$ be integers, 
and let $G$ be a non-complete graph. 
If $\sigma_{2}(G) \ge 2(s_{1} + s_{2} + 1) - 1$, 
then 
there exist two disjoint induced subgraphs $H_{1}$ and $H_{2}$ of $G$ 
such that 
$\sigma_{2}(H_{i}) \ge 2s_{i} - 1$ 
and 
$|H_{i}| \ge s_{i} + 1$ for $i \in \{ 1, 2 \}$. 
\end{thm}

\begin{thm} 
\label{disjoint subgraphs sigma2 tri-free}
Let $s_{1}, s_{2} \ge 2$ be integers, 
and let $G$ be a graph of order at least $3$. 
If $\sigma_{2}(G) \ge 2(s_{1} + s_{2}) - 1$ and $g(G) \ge 4$, 
then 
there exist two disjoint induced subgraphs $H_{1}$ and $H_{2}$ of $G$ 
such that 
$\sigma_{2}(H_{i}) \ge 2s_{i} - 1$ 
and 
$|H_{i}| \ge 2s_{i}$ for $i \in \{ 1, 2 \}$. 
\end{thm}

The graphs $G_{1}$ and $G_{2}$ defined on the previous page also show that 
the constraints on $\sigma_{2}$ in Theorems~\ref{disjoint subgraphs sigma2} and \ref{disjoint subgraphs sigma2 tri-free} cannot be weakened.

In order to show Theorems~\ref{disjoint subgraphs sigma2} and \ref{disjoint subgraphs sigma2 tri-free}, 
we actually prove slightly stronger results as follows. 
Here, for a graph $G$ and an integer $s$, 
we define $V_{\le s}(G) = \{v \in V(G) : d_{G}(v) \le s\}$.

\begin{thm} 
\label{feasible pair}
Let $s_{1}, s_{2} \ge 2$ be integers, 
and let $G$ be a non-complete graph. 
If $\sigma_{2}(G) \ge 2(s_{1} + s_{2} + 1) - 1$, 
then there exist two disjoint induced subgraphs $H_{1}$ and $H_{2}$ of $G$ 
such that for each $i$ with $i \in \{1, 2\}$, the following hold: 
\begin{enumerate}[{\upshape(i)}]
\setlength{\parskip}{1pt}
\setlength{\itemsep}{1pt}
\item $d_{H_{i}}(u) \ge s_{i}$ for $u \in V(H_{i}) \setminus V_{\le s_{1} + s_{2}}(G)$.  
\item $d_{H_{i}}(u) + d_{H_{i}}(v) \ge 2s_{i} - 1$ 
for $u \in V(H_{i}) \setminus V_{\le s_{1} + s_{2}}(G)$ 
and $v \in V(H_{i}) \cap V_{\le s_{1} + s_{2}}(G)$ 
with $uv \notin E(H_{i})$. 
\item $|H_{i}| \ge s_{i} + 1$. 
\end{enumerate}
\end{thm}

\begin{thm} 
\label{feasible pair tri.-free}
Let $s_{1}, s_{2} \ge 2$ be integers, 
and let $G$ be a graph of order at least $3$. 
If $\sigma_{2}(G) \ge 2(s_{1} + s_{2}) - 1$ and $g(G) \ge 4$, 
then there exist two disjoint induced subgraphs $H_{1}$ and $H_{2}$ of $G$ 
such that for each $i$ with $i \in \{1, 2\}$, the following hold: 
\begin{enumerate}[{\upshape(i)}]
\setlength{\parskip}{1pt}
\setlength{\itemsep}{1pt}
\item $d_{H_{i}}(u) \ge s_{i}$ for $u \in V(H_{i}) \setminus V_{\le s_{1} + s_{2}-1}(G)$.  
\item $d_{H_{i}}(u) + d_{H_{i}}(v) \ge 2s_{i} - 1$ 
for $u \in V(H_{i}) \setminus V_{\le s_{1} + s_{2}-1}(G)$ 
and $v \in V(H_{i}) \cap V_{\le s_{1} + s_{2}-1}(G)$ 
with $uv \notin E(H_{i})$. 
\item $|H_{i}| \ge 2s_{i}$. 
\end{enumerate}
\end{thm}

Note that 
if $G$ is a graph with $\sigma_{2}(G) \ge 2(s_{1} + s_{2} + 1) - 1$, 
then $G[V_{\le s_{1}+s_{2}}(G)]$ is a complete graph 
(see also Lemma~\ref{sigma}(\ref{complete}) in Subsection~\ref{Terminology and notation}). 
Therefore, 
for any two distinct non-adjacent vertices in such a graph $G$, 
at least one of the two vertices belongs to $V(G) \setminus V_{\le s_{1}+s_{2}}(G)$, 
i.e, 
(i) and (ii) of Theorem~\ref{feasible pair} imply that $\sigma_{2}(H_{i}) \ge 2s_{i} - 1$. 
Thus 
Theorem~\ref{disjoint subgraphs sigma2} 
immediately follows from Theorem~\ref{feasible pair}. 
Moreover, 
since $V_{\le s_{1}+s_{2}}(G) = \emptyset$ if and only if 
$\delta(G) \ge s_{1} + s_{2} + 1$ 
for a graph $G$, 
Theorem~\ref{feasible pair} also implies Step~1 for Theorem~\ref{Stiebitz}. 
Similarly, 
Theorem~\ref{feasible pair tri.-free} implies Theorem~\ref{disjoint subgraphs sigma2 tri-free} 
and Step~1 for Theorem~\ref{Kaneko}. 
In the next section, 
we give some concepts and lemmas to prove Theorems~\ref{feasible pair} and \ref{feasible pair tri.-free}. 
We will prove Theorems~\ref{feasible pair} and \ref{feasible pair tri.-free} 
in Sections~\ref{proof} and \ref{proof of tri.-free}, respectively.

These kinds of results are sometimes useful tools to get degree conditions for 
\textit{packing} of graphs, 
i.e., 
the existence of $k$ disjoint subgraphs which belong to some fixed class of graphs. 
In the last section (Section~\ref{applications}), 
we will explain it by taking disjoint cycles for example,  
and give some corollaries about it.

\medskip

We mention similar results. 
In 1966, Lov\'{a}sz \cite{Lovasz} 
proved a dual type of Theorem~\ref{Stiebitz} with respect to maximum degree: 
Every graph with maximum degree at most $s_{1} + s_{2} + 1$ 
has a partition $(H_{1}, H_{2})$ such that 
the maximum degree of each $H_{i}$ is at most $s_{i}$. 
On the other hand, 
Thomassen \cite{Thomassen2, Thomassen3} conjectured the connectivity version of Theorem~\ref{Stiebitz}: 
Every $(s_{1}+s_{2}+1)$-connected graph has 
a partition $(H_{1}, H_{2})$ such that 
each $H_{i}$ is $s_{i}$-connected, 
and 
he showed that this conjecture is true for $s_{2} \le 2$ (see \cite{Thomassen1}). 
However, this conjecture is still wide open for other cases, 
and hence there is a huge gap between ``degree'' and ``connectivity''. 
Other similar concepts can be found in \cite{Dunbar-Frick, Dunbar-Frick-Bullock, Kuhn-Osthus, Yang-Vumar}. 
Therefore, this type of problem has been extensively studied.

\section{Preparations for the proofs of Theorems~\ref{feasible pair} and \ref{feasible pair tri.-free}}
\label{preparation}

\subsection{Terminology and notation} 
\label{Terminology and notation}

We first prepare terminology and notation 
which we use in the rest of this paper, 
and we also give some lemmas in this subsection (Lemmas~\ref{sigma} and \ref{w(G_{1}, G_{2})}).

A \textit{clique} of a graph $G$ is a (possibly empty) vertex subset of $G$ that induces a complete subgraph of $G$, 
and we denote by $\omega(G)$ the cardinality of the largest clique of $G$. 
In this paper, 
we often identify a subgraph $H$ of a graph $G$ with its vertex set $V(H)$. 
For example, 
we write $G-H$ instead of $G-V(H)$.

By the definition of $\sigma_{2}(G)$, we can obtain the following. 
Since the proof is easy, we omit it.

\begin{lem} 
\label{sigma}
Let $s \ge 1$ be an integer, 
and let $G$ be a graph 
with $\sigma_{2}(G) \ge 2s - 1$. 
Then, 
\begin{enumerate}[{\upshape(i)}]
\item
\label{complete}
$V_{\le s - 1}(G)$ is a clique 
$($and hence, 
if $g(G) \ge 4$, then $|V_{\le s - 1}(G)| \le 2)$. 
\item
\label{order of G}
If $G$ is non-complete, then $|G| \ge s + 2$. 
\item
\label{order of G 2}
If $G$ is non-complete and $g(G) \ge 4$, then $|G| \ge 2s$. 
\end{enumerate}
\end{lem}

For a partition $(G_{1}, G_{2})$ of a graph $G$ and integers $s_{1}, s_{2} \ge 1$, 
we define 
\begin{align*}
f(G_{1}, G_{2}, s_{1}, s_{2}) = |E(G_{1})| + |E(G_{2})| + s_{2}|G_{1}| + s_{1}|G_{2}|. 
\end{align*}
Then, by the definition of the function $f$, 
we can obtain the following. 
(In order to find disjoint subgraphs as in Theorems~\ref{feasible pair} and \ref{feasible pair tri.-free}, 
we will consider some partition which was chosen 
so that value $f$ is maximized, 
see Sections~\ref{proof} and \ref{proof of tri.-free} for more details.)

\begin{lem} 
\label{w(G_{1}, G_{2})}
Let $s_{1}, s_{2} \ge 1$ be integers, 
and let $(G_{1}, G_{2})$ be a partition of a graph $G$. 
If $u_{1}$ is a vertex of $G_{1}$ such that $d_{G_{2}}(u_{1}) - d_{G_{1}}(u_{1}) \ge s_{2} - s_{1} + k$ for some integer $k$, 
then $f(G_{1}', G_{2}', s_{1}, s_{2}) \ge f(G_{1}, G_{2}, s_{1}, s_{2}) + k$, 
where $G_{1}' = G_{1} - u_{1}$ and $G_{2}' = G[V(G_{2}) \cup \{u_{1}\}]$. 
\end{lem}
\proof 
By the definitions of $f, G_{1}'$ and $G_{2}'$, 
we get 
\begin{align*}
f(G_{1}', G_{2}', s_{1}, s_{2}) 
&
= |E(G_{1}')| + |E(G_{2}')| + s_{2}|G_{1}'| + s_{1}|G_{2}'| \\
&
=\big( |E(G_{1})| - d_{G_{1}}(u_{1}) \big) 
+ \big( |E(G_{2})| + d_{G_{2}}(u_{1}) \big) 
+ s_{2} ( |G_{1}| - 1) + s_{1} (|G_{2}| + 1) \\
&
= f(G_{1}, G_{2}, s_{1}, s_{2}) + 
\big( d_{G_{2}}(u_{1}) - d_{G_{1}}(u_{1}) \big) - (s_{2} - s_{1}).  
\end{align*}
Hence, if 
$d_{G_{2}}(u_{1}) - d_{G_{1}}(u_{1}) \ge s_{2} - s_{1} + k$, 
then the assertion clearly holds. 
\qed

\subsection{Feasible graphs and degenerate graphs} 
\label{Feasible graphs and degenerate graphs}

In this subsection, 
we generalize the concepts of feasible graphs and degenerate graphs 
which were used in Stiebitz' argument \cite{Stiebitz}, 
and we will give some remarks and lemmas about it.

Now, let $G$ be a graph, and let $X$ be a clique of $G$. 
For an integer $s \ge 1$, 
an induced subgraph $H$ of $G$ is said to be \textit{$(s; X)$-feasible} in $G$, 
if $H$ satisfies the following conditions (F\ref{f1})--(F\ref{f3}): 
\begin{enumerate}
\renewcommand{\labelenumi}{\upshape{(F\arabic{enumi}})}
\item
\label{f1} 
$d_{H}(u) > s$ for $u \in V(H) \setminus X$. 
\item
\label{f2} 
$d_{H}(u) + d_{H}(v) > 2s$ for $u \in V(H) \setminus X$ and $v \in V(H) \cap X$ with $uv \notin E(H)$. 
\item
\label{f3} 
$H$ is non-complete.

(Note that the conditions (F\ref{f1}) and (F\ref{f2}) imply $\sigma_{2}(H) > 2s$ 
because $X$ is a clique.)
\end{enumerate}
For an integer $s \ge 1$, 
an induced subgraph $G'$ of $G$ is said to be \textit{$(s; X)$-degenerate} in $G$ 
if $G'$ satisfies the following conditions (D\ref{d1}) and (D\ref{d2}): 
\begin{enumerate}
\renewcommand{\labelenumi}{\upshape{(D\arabic{enumi}})}
\item
\label{d1} 
If $H$ is any non-complete induced subgraph of $G'$, 
then $H$ is not $(s; X)$-feasible in $G$, 
i.e., 
one of the following holds: 
\begin{itemize}
\item
$d_{H}(u) \le s$ for some $u \in V(H) \setminus X$, or
\item
$d_{H}(u) + d_{H}(v) \le 2s$ for some $u \in V(H) \setminus X$ and $v \in V(H) \cap X$ with $uv \notin E(H)$.
\end{itemize}
\item
\label{d2} 
$G'$ is non-complete. 
\end{enumerate}
For integers $s_{1}, s_{2} \ge 1$, 
a partition $(G_{1}, G_{2})$ of $G$ is called an $(s_{1}, s_{2}; X)$-\textit{degenerate partition} of $G$ 
if each $G_{i}$ is $(s_{i}; X)$-degenerate.

\begin{rem} 
\label{definition corollary}
Let $s \ge 1$ be an integer and $X$ be a clique in a graph $G$. 
Then, a non-complete induced subgraph $G'$ of $G$ is $(s; X)$-degenerate 
if and only if 
$G'$ contains no $(s; X)$-feasible induced subgraph. 
\end{rem}

\begin{lem} 
\label{feasible and degenerate}
Let $s \ge 2$ be an integer and $X$ be a clique in a graph $G$. 
Suppose that $\omega (G) \le s$. 
Then, the following hold. 
\begin{enumerate}[{\upshape(i)}]
\item
\label{(s; X)-feasible to (s-1; X)-feasible}
If 
$H$ is an $(s; X)$-feasible induced subgraph in $G$, 
then 
$H - x$ is an $(s-1; X)$-feasible induced subgraph in $G$ 
for every vertex $x$ of $H$. 
\item
\label{(s-1; X)-dege to (s; X)-dege}
If 
$G'$ is an $(s-1; X)$-degenerate induced subgraph in $G$, 
then 
$G[V(G') \cup \{x\}]$ is an $(s; X)$-degenerate induced subgraph in $G$ 
for every vertex $x$ of $G-G'$. 
\end{enumerate}
\end{lem}
\proof 
(i) Let $x \in V(H)$, and let $H' = H - x$. 
Then by (F\ref{f1}), 
$d_{H'}(u) > s - |\{x\}| = s - 1$ for $u \in V(H') \setminus X$. 
Similarly, by (F\ref{f2}), 
we also have 
$d_{H'}(u) + d_{H'}(v) > 2s - 2|\{x\}| = 2(s - 1)$ 
for 
$u \in V(H') \setminus X$ 
and 
$v \in V(H') \cap X$ 
with $uv \notin E(H')$. 
Furthermore, 
by (F\ref{f1})--(F\ref{f3}) and Lemma~\ref{sigma}(\ref{order of G}), 
$|H| \ge s + 2$, 
and hence $|H'| \ge s + 1$. 
This, together with $\omega(G) \le s$, 
implies that $H'$ is non-complete. 
Thus (\ref{(s; X)-feasible to (s-1; X)-feasible}) is proved.

(ii) To show (\ref{(s-1; X)-dege to (s; X)-dege}), 
let $x \in V(G-G')$, 
and let $G'' = G[V(G') \cup \{x\}]$. 
Note that by (D\ref{d2}), $G''$ is also non-complete. 
Suppose that $G''$ is not $(s; X)$-degenerate. 
Then by Remark~\ref{definition corollary}, 
$G''$ contains an $(s; X)$-feasible induced subgraph $H$. 
If $V(H) \subseteq V(G')$, 
then this contradicts that $G'$ is $(s-1; X)$-degenerate. 
Thus $x \in V(H)$. 
But then, by Lemma~\ref{feasible and degenerate}(\ref{(s; X)-feasible to (s-1; X)-feasible}), 
$H-x$ is an $(s-1; X)$-feasible induced subgraph in $G'$, 
a contradiction again. 
Thus (\ref{(s-1; X)-dege to (s; X)-dege}) is proved. 
\qed

\begin{lem} 
\label{s_{i} + l + 1}
Let $s_{1}, s_{2} \ge 1$ be integers and $\varepsilon \in \{0, 1\}$, 
and let $G$ be a graph with $\sigma_{2}(G) \ge 2(s_{1} + s_{2} + \varepsilon) - 1$. 
Furthermore, 
let $X = V_{\le s_{1} + s_{2} + \varepsilon-1}(G)$ 
and suppose $(G_{1}, G_{2})$ is a partition of $G$ 
such that $G_{1}$ is $(s_{1} - 1 + \varepsilon; X)$-degenerate. 
Then, there is a vertex $u_{1}$ of $G_{1}$ such that 
$d_{G_{2}}(u_{1}) - d_{G_{1}}(u_{1}) \ge s_{2} - s_{1} + 1$. 
\end{lem}
\proof 
We consider two cases. 
\begin{enumerate}[{\textup{{\bf Case~\arabic{enumi}.}}}]
\setcounter{enumi}{0}
\item $d_{G_{1}}(u) \le s_{1} + \varepsilon - 1$ for some $u \in V(G_{1}) \setminus X$.  
\end{enumerate}

Since $u \notin X$, we have $d_{G}(u) \ge s_{1} + s_{2} + \varepsilon$, 
and hence $d_{G_{2}}(u) \ge s_{2} + 1$. 
Thus, we get 
\begin{align*}
d_{G_{2}}(u) - d_{G_{1}}(u) 
\ge (s_{2} + 1) - (s_{1} + \varepsilon - 1) 
= s_{2} - s_{1} + 2 - \varepsilon 
\ge s_{2} - s_{1} + 1, 
\end{align*}
so in this case, $u$ is the desired vertex $u_{1}$.

\begin{enumerate}[{\textup{{\bf Case~\arabic{enumi}.}}}]
\setcounter{enumi}{1}
\item $d_{G_{1}}(u) \ge s_{1} + \varepsilon$ for every $u \in V(G_{1}) \setminus X$.  
\end{enumerate}

Then, since $G_{1}$ is $(s_{1} - 1 + \varepsilon; X)$-degenerate, 
there exist two non-adjacent vertices $u \in V(G_{1}) \setminus X$ and $v \in V(G_{1}) \cap X$ 
such that $d_{G_{1}}(u) + d_{G_{1}}(v) \le 2 (s_{1} - 1 + \varepsilon)$. 
By the assumption of Case~2, 
we have $d_{G_{1}}(v) \le s_{1} -2 + \varepsilon$, 
so $d_{G_{1}}(v) = s_{1} -l + \varepsilon$ for some integer $l \ge 2$.

\begin{enumerate}[{\textup{{\bf Case~2.\arabic{enumi}.}}}]
\setcounter{enumi}{0}
\item $d_{G_{2}}(v) \ge s_{2} - l + \varepsilon + 1$.  
\end{enumerate}

In this case, 
\begin{align*}
d_{G_{2}}(v) - d_{G_{1}}(v) \ge (s_{2} - l + \varepsilon + 1) - (s_{1} - l + \varepsilon) = s_{2} - s_{1} + 1, 
\end{align*}
so $v$ is the desired vertex $u_{1}$.

\begin{enumerate}[{\textup{{\bf Case~2.\arabic{enumi}.}}}]
\setcounter{enumi}{1}
\item $d_{G_{2}}(v) \le s_{2} - l + \varepsilon$.  
\end{enumerate}

Then, 
\begin{align*}
d_{G}(v) = d_{G_{1}}(v) + d_{G_{2}}(v) \le s_{1} + s_{2} - 2l + 2\varepsilon \ (\le s_{1} + s_{2} - 2). 
\end{align*}
Note that 
$d_{G_{1}}(u) \le s_{1} - 2 + l + \varepsilon$ because 
$d_{G_{1}}(u) + d_{G_{1}}(v) \le 2(s_{1} - 1 + \varepsilon)$ and $d_{G_{1}}(v) = s_{1} - l + \varepsilon$, 
and hence 
\begin{align*}
d_{G_{2}}(u) 
&\ge \sigma_{2}(G) - d_{G_{1}}(u) - d_{G}(v) \\
&\ge \big( 2(s_{1} + s_{2} + \varepsilon) - 1 \big) - (s_{1} - 2 + l + \varepsilon) - (s_{1} + s_{2} - 2l + 2\varepsilon)\\
&= s_{2} + l + 1 - \varepsilon = (s_{2} - 2 + l +\varepsilon) + 3 - 2 \varepsilon \ge (s_{2} - 2 + l + \varepsilon) + 1. 
\end{align*}
This implies that $d_{G_{2}}(u) - d_{G_{1}}(u) \ge s_{2} - s_{1} + 1$. 
Thus, $u$ is the desired vertex $u_{1}$. 
\qed

\subsection{Lemmas for Theorem~\ref{feasible pair}} 
\label{lemmas for feasible pair}

In this subsection, 
we prepare lemmas which will be used in the proof of Theorem~\ref{feasible pair}.

\begin{lem} 
\label{|U| = s_{i}+alpha}
Let $s_{1}, s_{2} \ge 2$ be integers 
and $G$ be a non-complete graph with $\sigma_{2}(G) \ge 2 (s_{1} + s_{2} + 1) - 1$, 
and let $X = V_{\le s_{1} + s_{2}}(G)$. 
If $U$ is a vertex subset of $G$ such that $|U| = s_{2} + \alpha$ for some $\alpha \le 1$, 
then the following hold: 
{\rm (i)} $d_{G-U}(u) > s_{1}-\alpha$ for $u \in V(G-U) \setminus X$, 
{\rm (ii)} $d_{G-U}(u) + d_{G-U}(v) > 2(s_{1}-\alpha)$ 
for $u \in V(G-U) \setminus X$ 
and 
$v \in V(G-U) \cap X$ with $uv \notin E(G-U)$, 
and 
{\rm (iii)} $|G-U| \ge (s_{1}-\alpha) + 3$. 
\end{lem}
\proof 
By the degree condition of $G$, 
the definition of $X$ and the assumption that $|U| = s_{2} + \alpha$, 
we can easily check that (i) and (ii) hold. 
Moreover, 
since $G$ is non-complete, 
it follows from Lemma~\ref{sigma}(\ref{order of G}) that 
$|G| \ge s_{1} + s_{2} + 3$, 
and hence $|G-U| \ge (s_{1}-\alpha) + 3$. 
Thus (iii) also holds. 
\qed

\begin{lem} 
\label{G-G' is non-complete}
Let $s_{1}, s_{2} \ge 2$ be integers 
and $G$ be a non-complete graph 
with $\sigma_{2}(G) \ge 2 (s_{1} + s_{2} + 1) - 1$ 
and 
$\omega(G) \le \min \{s_{1}, s_{2}\}$, 
and let $X = V_{\le s_{1} + s_{2}}(G)$. 
If $G'$ is an $(s_{1}; X)$-degenerate induced subgraph in $G$, 
then $G-G'$ is non-complete. 
\end{lem}
\proof
Note that $G'$ is a proper subgraph of $G$ 
because $G$ is not $(s_{1}; X)$-degenerate. 
Suppose that $G-G'$ is complete.  
Since $|G-G'| \le \omega (G) \le s_{2}$, 
we have $|G-G'| = s_{2} + \alpha$ for some integer $\alpha \le 0$. 
Note that $(s_{1}-\alpha) + 3 > s_{1} \ge \omega(G)$. 
By applying Lemma~\ref{|U| = s_{i}+alpha} with $U = V(G-G')$, 
it follows that 
$G' = G-U$ is an $(s_{1} - \alpha)$-feasible induced subgraph in $G$, 
which contradicts that $G'$ is  $(s_{1}; X)$-degenerate. 
\qed

\subsection{Lemmas for Theorem~\ref{feasible pair tri.-free}} 
\label{lemmas for feasible pair tri.-free}

In this subsection, 
we prepare lemmas which will be used in the proof of Theorem~\ref{feasible pair tri.-free}.

\begin{lem} 
\label{(s-alpha)-feasible subgraph}
Let $s_{1}, s_{2} \ge 2$ be integers 
and $G$ be a graph of order at least $3$ 
with $\sigma_{2}(G) \ge 2 (s_{1} + s_{2}) - 1$ and $g(G) \ge 4$, 
and let $X = V_{\le s_{1} + s_{2} - 1}(G)$. 
If $U$ is a vertex subset of $G$ such that 
$|U| \le 3$ 
and 
$d_{G[U]}(x) \le 1 + \alpha$ for some $\alpha \in \{0, 1\}$ and for all $x \in V(G) \setminus U$, 
then 
$G-U$ is an $(s_{1}-\alpha; X)$-feasible induced subgraph in $G$. 
\end{lem}
\proof 
By the degree condition of $G$ and the definition of $X$, 
$d_{G-U}(u) 
\ge 
d_{G}(u) - d_{G[U]}(u) 
\ge (s_{1} + s_{2}) - (1 + \alpha) 
= (s_{1} - \alpha) + s_{2} - 1 > s_{1} - \alpha$ for $u \in V(G-U) \setminus X$. 
Similarly, 
we have 
$d_{G-U}(u) + d_{G-U}(v) 
\ge \sigma_{2}(G) - \big( d_{G[U]}(u) + d_{G[U]}(v) \big) 
\ge \big( 2(s_{1} + s_{2})-1 \big) - 2(1 + \alpha) 
= 2(s_{1} - \alpha) + 2s_{2} - 3 > 2(s_{1} - \alpha)$ for 
$u \in V(G-U) \setminus X$ and $v \in V(G-U) \cap X$ with $uv \notin E(G-U)$. 
By Lemma~\ref{sigma}(\ref{order of G 2}), $|G| \ge 2(s_{1} + s_{2}) > 6$, 
and hence $|G-U| \ge 3$. 
Since $g(G) \ge 4$, this implies that $G-U$ is non-complete. 
Thus, $G-U$ is $(s_{1}-\alpha; X)$-feasible in $G$. 
\qed

\begin{lem} 
\label{G-G' is non-complete tri.-free}
Let $s_{1}, s_{2} \ge 2$ be integers 
and $G$ be a graph of order at least $3$ 
with $\sigma_{2}(G) \ge 2 (s_{1} + s_{2}) - 1$ and $g(G) \ge 4$, 
and let $X = V_{\le s_{1} + s_{2} - 1}(G)$. 
If $G'$ is an $(s_{1}; X)$-degenerate induced subgraph in $G$, 
then 
$G-G'$ is non-complete. 
\end{lem}
\proof 
Note that $G'$ is a proper subgraph of $G$  
because $G$ itself is not $(s_{1}; X)$-degenerate. 
Suppose that $G-G'$ is complete. 
Then by Lemma~\ref{sigma}(\ref{complete}), 
$|G-G'| \le 2$. 
Since $g(G) \ge 4$, 
this implies that 
$d_{G-G'}(x) \le 1$ 
for all $x \in V(G') = V(G) \setminus V(G-G')$. 
Therefore, 
by applying Lemma~\ref{(s-alpha)-feasible subgraph} 
with $U = V(G-G')$ and $\alpha = 0$, 
it follows that $G' = G-U$ is an $(s_{1}; X)$-feasible induced subgraph, a contradiction. 
\qed

\begin{lem} 
\label{(s-1)-feasible subgraph}
Let $s_{1}, s_{2} \ge 2$ be integers 
and $G$ be a graph of order at least $3$ 
with $\sigma_{2}(G) \ge 2 (s_{1} + s_{2}) - 1$ and $g(G) \ge 4$, 
and let $X = V_{\le s_{1} + s_{2} - 1}(G)$. 
Furthermore, 
let $(G_{1}, G_{2})$ be an $(s_{1}, s_{2}-1; X)$-degenerate partition 
or 
an $(s_{1}-1, s_{2}; X)$-degenerate partition of $G$ 
such that $f(G_{1}, G_{2}, s_{1}, s_{2})$ is maximum. 
If $G_{1}$ is an $(s_{1}-1; X)$-degenerate induced subgraph in $G$, 
then 
$G_{2}$ is not an $(s_{2}-1; X)$-degenerate induced subgraph in $G$. 
\end{lem}
\proof 
Suppose that $G_{1}$ is an $(s_{1}-1; X)$-degenerate induced subgraph in $G$. 
Note that by the assumption of Lemma~\ref{(s-1)-feasible subgraph}, 
$G_{2}$ is $(s_{2}; X)$-degenerate. 
By applying Lemma~\ref{s_{i} + l + 1} with $\varepsilon = 0$, 
we can take a vertex $u_{1}$ of $G_{1}$ 
such that $d_{G_{2}}(u_{1}) - d_{G_{1}}(u_{1}) \ge s_{2} - s_{1} + 1$. 
Then by Lemma~\ref{w(G_{1}, G_{2})}, 
$f(G_{1}', G_{2}', s_{1}, s_{2}) > f(G_{1}, G_{2}, s_{1}, s_{2})$, 
where $G_{1}' = G_{1} - u_{1}$ and $G_{2}' = G[V(G_{2}) \cup \{u_{1}\}]$. 
Note that $G_{2}'$ is non-complete because $G_{2}$ is non-complete.

Assume first that $G_{1}'$ is non-complete. 
Then, since $G_{1}$ is $(s_{1}-1; X)$-degenerate, 
$G_{1}'$ is clearly $(s_{1}-1; X)$-degenerate. 
Hence by the choice of $(G_{1}, G_{2})$, 
$G_{2}'$ is not $(s_{2}; X)$-degenerate. 
Therefore, by Remark~\ref{definition corollary}, 
$G_{2}'$ contains an $(s_{2}; X)$-feasible induced subgraph $H$. 
If $u_{1} \notin V(H)$, 
then $H$ is a subgraph of $G_{2}$, which contradicts that $G_{2}$ is $(s_{2}; X)$-degenerate. 
Thus, we have $u_{1} \in V(H)$. 
Note that $\omega(G) \le 2 \le s_{2}$ because $g(G) \ge 4$. 
Hence, by Lemma~\ref{feasible and degenerate}(\ref{(s; X)-feasible to (s-1; X)-feasible}), 
$H-u_{1}$ is an induced subgraph of $G_{2}$ which is $(s_{2}-1; X)$-feasible. 
Therefore, by Remark~\ref{definition corollary}, 
$G_{2}$ is not $(s_{2}-1; X)$-degenerate.

Assume next that $G_{1}'$ is complete. 
By Lemma~\ref{sigma}(\ref{complete}), $|G_{1}'| \le 2$, 
and hence we have $|G_{1}| \le 3$. 
Since $g(G) \ge 4$, this implies that $d_{G_{1}}(x) \le 2$ for all $x \in V(G_{2})$. 
Therefore, by applying Lemma~\ref{(s-alpha)-feasible subgraph} with $U = V(G_{1})$ and $\alpha = 1$, 
and by using the symmetry of $s_{1}, s_{2}$ in Lemma~\ref{(s-alpha)-feasible subgraph}, 
it follows that $G_{2} = G - U$ is $(s_{2}-1; X)$-feasible, 
that is, $G_{2}$ is not $(s_{2}-1; X)$-degenerate.
\qed

\section{Proof of Theorem~\ref{feasible pair}} 
\label{proof}

Let $s_{1}, s_{2} \ge 2$ be integers, 
and let $G$ be a non-complete graph with $\sigma_{2}(G) \ge 2 (s_{1} + s_{2} + 1) - 1$. 

\paragraph{}
First assume that $\omega(G) > \min \{s_{1}, s_{2}\}$. 
By the symmetry of $s_{1}$ and $s_{2}$, we may assume that 
$G$ contains a clique $U$ of order $s_{2} + 1$. 
Then, by Lemma~\ref{|U| = s_{i}+alpha}, 
we can easily see that 
$H_{1} = G - U$ and $H_{2} = G[U]$ are desired subgraphs.

\paragraph{}
Now assume that $\omega(G) \le \min \{s_{1}, s_{2}\}$. 
Let $X = V_{\le s_{1} + s_{2}}(G)$. 
To finish the proof, it suffices to find a partition $(G_{1}^{*}, G_{2}^{*})$ of $G$ 
such that 
\begin{align*}
\textup{each $G_{i}^{*}$ contains an $(s_{i}-1; X)$-feasible induced subgraph.} 
\end{align*}
To find such a partition, we distinguish two cases.

\begin{enumerate}[{\textup{{\bf Case~\arabic{enumi}.}}}]
\setcounter{enumi}{0}
\item $G$ has neither an $(s_{1}, s_{2}-1; X)$-degenerate partition nor an $(s_{1}-1, s_{2}; X)$-degenerate partition.  
\end{enumerate}

Let $G_{1}$ be an $(s_{1}-1; X)$-feasible induced subgraph of smallest possible order in $G$. 
Then, it follows from Lemma~\ref{feasible and degenerate}(\ref{(s; X)-feasible to (s-1; X)-feasible}) 
that $G_{1}$ does not contain any $(s_{1}; X)$-feasible induced subgraph, 
i.e., $G_{1}$ is an $(s_{1}; X)$-degenerate (by Remark~\ref{definition corollary}). 
Hence, by the assumption of Case 1, 
$G_{2} := G - G_{1}$ is not $(s_{2}-1; X)$-degenerate. 
Note that, 
by Lemma~\ref{G-G' is non-complete}, 
$G_{2} \ (= G - G_{1})$ is non-complete, 
and hence 
it follows from Remark~\ref{definition corollary} that 
$G_{2}$ contains an $(s_{2}-1; X)$-feasible induced subgraph. 
Thus, $(G_{1}, G_{2})$ is the desired partition $(G_{1}^{*}, G_{2}^{*})$.

\begin{enumerate}[{\textup{{\bf Case~\arabic{enumi}.}}}]
\setcounter{enumi}{1}
\item $G$ contains an $(s_{1}, s_{2}-1; X)$-degenerate partition or 
an $(s_{1}-1, s_{2}; X)$-degenerate partition.  
\end{enumerate}

Among all such partitions, 
let $(G_{1}, G_{2})$ be one such that 
$f(G_{1}, G_{2}, s_{1}, s_{2})$ is as large as possible 
(see the definition of $f$ in Subsection~\ref{Terminology and notation}). 
By symmetry, we may assume that 
$(G_{1}, G_{2})$ is an $(s_{1}, s_{2}-1; X)$-degenerate partition of $G$. 
We now apply Lemma~\ref{s_{i} + l + 1} with $\varepsilon = 1$. 
Then, we can take a vertex $u_{1}$ of $G_{1}$ 
such that 
$d_{G_{2}}(u_{1}) - d_{G_{1}}(u_{1}) \ge s_{2} - s_{1} + 1$, 
and hence by Lemma~\ref{w(G_{1}, G_{2})}, 
$f(G_{1}', G_{2}', s_{1}, s_{2}) > f(G_{1}, G_{2}, s_{1}, s_{2})$, 
where $G_{1}' = G_{1} - u_{1}$ and $G_{2}' = G[V(G_{2}) \cup \{u_{1}\}]$. 
Note that $G_{2}'$ is non-complete. 
In the rest of this proof, 
we show that $(G_{1}', G_{2}')$ is the desired partition $(G_{1}^{*}, G_{2}^{*})$.

Since $G_{2}' = G[V(G_{2}) \cup \{u_{1}\}]$ 
and $G_{2}$ is $(s_{2}-1; X)$-degenerate, 
it follows from Lemma~\ref{feasible and degenerate}(\ref{(s-1; X)-dege to (s; X)-dege}) 
that 
$G_{2}'$ is $(s_{2}; X)$-degenerate. 
Then, by Lemma~\ref{G-G' is non-complete} and the symmetry of $s_{1}, s_{2}$ in Lemma~\ref{G-G' is non-complete}, 
$G_{1}' \ (= G - G_{2}')$ is non-complete. 
Hence, the choice of $(G_{1}, G_{2})$ and Remark~\ref{definition corollary} yield that 
$G_{1}'$ contains an $(s_{1}-1; X)$-feasible induced subgraph.

We next show that $G_{2}'$ contains an $(s_{2}-1; X)$-feasible induced subgraph. 
To show it, suppose that 
$G_{2}'$ is $(s_{2}-1; X)$-degenerate. 
Recall that $G_{1}'$ is non-complete as mentioned in the above. 
Then $G_{1}'$ is clearly $(s_{1}; X)$-degenerate 
because 
$G_{1}$ is $(s_{1}; X)$-degenerate. 
Thus, 
$(G_{1}', G_{2}')$ is an $(s_{1}, s_{2}-1; X)$-degenerate partition 
such that $f(G_{1}', G_{2}', s_{1}, s_{2}) > f(G_{1}, G_{2}, s_{1}, s_{2})$, 
which contradicts the choice of $(G_{1}, G_{2})$. 
Thus, $G_{2}'$ is not $(s_{2}-1; X)$-degenerate, 
and hence by Remark~\ref{definition corollary}, $G_{2}'$ contains an $(s_{2}-1; X)$-feasible induced subgraph. 
Therefore, $(G_{1}', G_{2}')$ is the desired partition $(G_{1}^{*}, G_{2}^{*})$.

This completes the proof of Theorem~\ref{feasible pair}. 
\qed

\section{Proof of Theorem~\ref{feasible pair tri.-free}} 
\label{proof of tri.-free}

Let $s_{1}, s_{2} \ge 2$ be integers, 
and let $G$ be a graph of order at least $3$ with $\sigma_{2}(G) \ge 2 (s_{1} + s_{2} + 1) - 1$ and $g(G) \ge 4$.

\paragraph{}
First assume that $|V_{\le s_{1} + s_{2} - 2}(G)| \ge 2$, 
say $v_{1}, v_{2} \in V_{\le s_{1} + s_{2} - 2}(G)$ with $v_{1} \neq v_{2}$. 
Note that by Lemma~\ref{sigma}(\ref{complete}), $v_{1}v_{2} \in E(G)$. 
Since $g(G) \ge 4$, 
this implies that $d_{G[\{v_{1}, v_{2}\}]}(u) \le 1$ for all $u \in V(G) \setminus \{v_{1}, v_{2}\}$. 
Let $u$ be an arbitrary vertex of $G - \{v_{1}, v_{2}\}$, 
and then we may assume that $uv_{1} \notin E(G)$. 
Hence, 
\begin{align*}
d_{G - \{v_{1}, v_{2}\}}(u) 
\ge \sigma_{2}(G) - d_{G}(v_{1}) - d_{G[\{v_{1}, v_{2}\}]}(u) 
\ge \big( 2(s_{1} + s_{2}) - 1\big) - (s_{1} + s_{2} - 2) - 1 = s_{1} + s_{2}. 
\end{align*}
Since $u$ is an arbitrary vertex of $G - \{v_{1}, v_{2}\}$, 
this implies that $\delta(G - \{v_{1}, v_{2}\}) \ge s_{1} + s_{2}$. 
Hence by Theorem~\ref{Kaneko}, 
there exists a partition $(H_{1}, H_{2})$ of $G- \{v_{1}, v_{2}\}$ 
such that $\delta(H_{i}) \ge s_{i}$ for $i \in \{1, 2\}$. 
This, together with Lemma~\ref{sigma}(\ref{order of G 2}), 
implies that 
$H_{1}$ and $H_{2}$ are desired induced subgraphs.

\paragraph{}
Now assume that 
\begin{align}
\label{degree of X}
|V_{\le s_{1} + s_{2} - 2}(G)| \le 1. 
\end{align}
Let 
$X = V_{\le s_{1} + s_{2}-1}(G)$. 
To finish the proof of Theorem~\ref{feasible pair tri.-free}, it suffices to find a partition $(G_{1}^{*}, G_{2}^{*})$ of $G$ 
such that 
\begin{align*}
\textup{each $G_{i}^{*}$ contains an $(s_{i}-1; X)$-feasible induced subgraph.} 
\end{align*}
To find such a partition, we distinguish two cases.

\begin{enumerate}[{\textup{{\bf Case~\arabic{enumi}.}}}]
\setcounter{enumi}{0}
\item $G$ has neither an $(s_{1}, s_{2}-1; X)$-degenerate partition nor an $(s_{1}-1, s_{2}; X)$-degenerate partition.  
\end{enumerate}

In this case, 
replacing ``Lemma~\ref{G-G' is non-complete}'' with ``Lemma~\ref{G-G' is non-complete tri.-free}'' 
in Case~1 of the proof of Theorem~\ref{feasible pair}, 
the same argument can work 
(note that we can use Lemma~\ref{feasible and degenerate}(\ref{(s; X)-feasible to (s-1; X)-feasible}) as $s = s_{1}$
because $\omega(G) \le 2 \le s_{1}$). 
Thus, we can obtain the desired conclusion.

\begin{enumerate}[{\textup{{\bf Case~\arabic{enumi}.}}}]
\setcounter{enumi}{1}
\item $G$ contains an $(s_{1}, s_{2}-1; X)$-degenerate partition or 
an $(s_{1}-1, s_{2}; X)$-degenerate partition.  
\end{enumerate}

Among all such partitions, 
let $(G_{1}, G_{2})$ be one such that 
\begin{enumerate}[]
\item{(C1)}
$f(G_{1}, G_{2}, s_{1}, s_{2})$ 
is as large as possible.  
\end{enumerate}
By symmetry, we may assume that 
$(G_{1}, G_{2})$ is an $(s_{1}, s_{2}-1; X)$-degenerate partition of $G$. 
Then, by Lemma~\ref{(s-1)-feasible subgraph} and the symmetry of $s_{1}, s_{2}$ in Lemma~\ref{(s-1)-feasible subgraph}, 
$G_{1}$ is not $(s_{1}-1; X)$-degenerate 
(and thus, $G_{1}$ and $G_{2}$ are not symmetric). 
We further choose such a partition so that 
\begin{enumerate}[]
\item{(C2)}
$|G_{1}|$ is as small as possible, subject to (C1).  
\end{enumerate}

We divide Case 2 into two parts.

\begin{enumerate}[{\textup{{\bf Case~2.\arabic{enumi}.}}}]
\setcounter{enumi}{0}
\item 
There exists a vertex $x$ of $G_{1}$ 
such that 
(i) $d_{G_{2}}(x) - d_{G_{1}}(x) \ge s_{2} - s_{1}$, 
and 
(ii) $G_{1} - x$ is not $(s_{1}-1; X)$-degenerate if $d_{G_{2}}(x) - d_{G_{1}}(x) = s_{2} - s_{1}$. 
\end{enumerate}

Let $d_{G_{2}}(x) - d_{G_{1}}(x) = s_{2} - s_{1} + k$ for some integer $k \ge 0$. 
Then, by Lemma~\ref{w(G_{1}, G_{2})}, 
$f(G_{1}', G_{2}', s_{1}, s_{2})  \ge f(G_{1}, G_{2}, s_{1}, s_{2}) + k$, 
where $G_{1}' = G_{1} - x$ and $G_{2}' = G[V(G_{2}) \cup \{x\}]$. 
Note that $G_{2}'$ is non-complete. 
We show that $(G_{1}', G_{2}')$ is the desired partition $(G_{1}^{*}, G_{2}^{*})$.

Since $\omega(G) \le 2 \le s_{2}$ 
and $G_{2}$ is $(s_{2}-1; X)$-degenerate, 
it follows from Lemma~\ref{feasible and degenerate}(\ref{(s-1; X)-dege to (s; X)-dege}) 
that $G_{2}'$ is $(s_{2}; X)$-degenerate. 
Then, by Lemma~\ref{G-G' is non-complete tri.-free} 
and the symmetry of $s_{1}, s_{2}$ in Lemma~\ref{G-G' is non-complete tri.-free}, 
$G_{1}' \ (= G - G_{2}')$ is non-complete. 
Therefore, 
if $k \ge 1$, then by the choice (C1) and Remark~\ref{definition corollary} 
yield that $G_{1}'$ contains an $(s_{1}-1; X)$-feasible induced subgraph; 
if $k = 0$, then  by the assumption of Case 2.1(ii) and Remark~\ref{definition corollary} 
yield that $G_{1}'$ contains an $(s_{1}-1; X)$-feasible induced subgraph. 
In either case, $G_{1}'$ contains an $(s_{1}-1; X)$-feasible induced subgraph.

We next show that $G_{2}'$ contains an $(s_{2}-1; X)$-feasible induced subgraph. 
To show it, suppose that 
$G_{2}'$ is $(s_{2}-1; X)$-degenerate. 
Recall that $G_{1}'$ is non-complete as mentioned in the above. 
Then, $G_{1}'$ is clearly $(s_{1}; X)$-degenerate 
because 
$G_{1}$ is $(s_{1}; X)$-degenerate. 
Thus, 
$(G_{1}', G_{2}')$ is an $(s_{1}, s_{2}-1; X)$-degenerate partition 
such that $f(G_{1}', G_{2}', s_{1}, s_{2}) \ge f(G_{1}, G_{2}, s_{1}, s_{2}) + k$ 
and $|G_{1}'| < |G_{1}|$, 
which contradicts (C1) or (C2). 
Thus, $G_{2}'$ is not $(s_{2}-1; X)$-degenerate, 
and hence by Remark~\ref{definition corollary}, $G_{2}'$ contains an $(s_{2}-1; X)$-feasible induced subgraph. 
Therefore, $(G_{1}', G_{2}')$ is the desired partition $(G_{1}^{*}, G_{2}^{*})$.

\begin{enumerate}[{\textup{{\bf Case~2.\arabic{enumi}.}}}]
\setcounter{enumi}{1}
\item 
Every vertex $x$ of $G_{1}$ 
satisfies that 
(i) $d_{G_{2}}(x) - d_{G_{1}}(x) \le s_{2} - s_{1}$, 
and 
(ii) $G_{1} - x$ is $(s_{1}-1; X)$-degenerate if $d_{G_{2}}(x) - d_{G_{1}}(x) = s_{2} - s_{1}$. 
\end{enumerate}

We define 
\begin{align*}
Z_{1} = \big\{ z \in V(G_{1}) \setminus X : d_{G_{1}}(z) = s_{1} \textup{ and } d_{G_{2}}(z) = s_{2} \big\}. 
\end{align*}

\begin{claim}
\label{Z is non-empty} 
We have $Z_{1} \neq \emptyset$. 
\end{claim}
\proof 
On the contrary, 
suppose that 
$Z_{1} = \emptyset$.

\begin{subclaim}
\label{sub1} 
$d_{G_{1}}(u) > s_{1}$ for $u \in V(G_{1}) \setminus X$. 
\end{subclaim}
\proof 
If there exists a vertex $u$ of $V(G_{1}) \setminus X$ 
such that $d_{G_{1}}(u) \le s_{1}$, 
then 
by the assumption of Case~2.2(i) 
and the definition of $X$, 
it is easy to check that 
$d_{G_{1}}(u) = s_{1}$ 
and 
$d_{G_{2}}(u) = s_{2}$, 
and hence $u \in Z_{1}$, a contradiction.  
\qed

Since $G_{1}$ is $(s_{1}; X)$-degenerate, 
it follows from Subclaim~\ref{sub1} that 
there exist two vertices 
$u_{1} \in V(G_{1}) \setminus X$ and $v_{1} \in V(G_{1}) \cap X$ 
with $u_{1}v_{1} \notin E(G_{1})$ 
such that $d_{G_{1}}(u_{1}) + d_{G_{1}}(v_{1}) \le 2s_{1}$. 
Since $d_{G_{1}}(u_{1}) > s_{1}$ by Subclaim~\ref{sub1}, 
we have $d_{G_{1}}(v_{1}) = s_{1} - l$ for some integer $l \ge 1$.

\begin{subclaim}
\label{sub2} 
$d_{G}(v_{1}) \le s_{1} +s_{2} - 2l - 1 \ ( \le s_{1} + s_{2} - 3)$. 
\end{subclaim}
\proof 
By the assumption of Case~2.2(i), we have $d_{G_{2}}(v_{1}) \le s_{2} - l$. 
Suppose that $d_{G_{2}}(v_{1}) = s_{2} - l$. 
Then 
$d_{G_{2}}(v_{1}) - d_{G_{1}}(v_{1}) = (s_{2} - l) - (s_{1} - l) = s_{2} - s_{1}$. 
Hence 
by the assumption Case~2.2(ii), 
$G_{1} - v_{1}$ is $(s_{1}-1; X)$-degenerate. 
Note that by Subclaim~\ref{sub1}, 
$d_{G_{1}-v_{1}}(u) > s_{1} - 1$ for $u \in V(G_{1} - v_{1}) \setminus X$. 
Therefore, there exist two vertices 
$u_{1}' \in V(G_{1} - v_{1}) \setminus X$ and $v_{1}' \in V(G_{1} - v_{1}) \cap X$ 
with $u_{1}'v_{1}' \notin E(G_{1} - v_{1})$ 
such that $d_{G_{1} - v_{1}}(u_{1}') + d_{G_{1} - v_{1}}(v_{1}') \le 2(s_{1}-1)$. 
Since $d_{G_{1}-v_{1}}(u_{1}') > s_{1} - 1$, 
this implies that $d_{G_{1} - v_{1}}(v_{1}') \le s_{1} - 2$, 
and hence $d_{G_{1}}(v_{1}') \le s_{1} - 1$. 
Combining this 
with the assumption of Case~2.2(i), 
we also get $d_{G_{2}}(v_{1}') \le s_{2} - 1$. 
Therefore, $v_{1}$ and $v_{1}'$ are two distinct vertices 
such that 
$d_{G}(v_{1}) = d_{G_{1}}(v_{1}) + d_{G_{2}}(v_{1}) 
\le (s_{1} - l) + (s_{2} - l) \le s_{1} + s_{2} - 2$ 
and 
$d_{G}(v_{1}') = d_{G_{1}}(v_{1}') + d_{G_{2}}(v_{1}') 
\le (s_{1} - 1) + (s_{2} - 1) =  s_{1} + s_{2} - 2$, 
which contradicts (\ref{degree of X}). 
Thus $d_{G_{2}}(v_{1}) \le s_{2} - l - 1$, 
and hence 
$d_{G}(v_{1}) = d_{G_{1}}(v_{1}) + d_{G_{2}}(v_{1}) 
\le s_{1} + s_{2}  - 2l - 1 \ ( \le s_{1} + s_{2}- 3)$. 
\qed

By Subclaim~\ref{sub2}, 
$d_{G}(u_{1}) \ge \sigma_{2}(G) - d_{G}(v_{1}) 
\ge \big( 2(s_{1} + s_{2}) - 1 \big) - (s_{1} + s_{2} - 2l - 1) = s_{1} + s_{2} + 2l$. 
On the other hand, 
since 
$d_{G_{1}}(u_{1}) + d_{G_{1}}(v_{1}) \le 2s_{1}$ and $d_{G_{1}}(v_{1}) = s_{1} - l$, 
it follows that 
$d_{G_{1}}(u_{1}) \le s_{1} + l$. 
Therefore, 
we get 
\begin{align*}
d_{G_{2}}(u_{1}) - d_{G_{1}}(u_{1}) 
&= \big( d_{G}(u_{1}) - d_{G_{1}}(u_{1})\big) - d_{G_{1}}(u_{1}) \\
&\ge \big( s_{1} + s_{2} + 2l - (s_{1} + l)\big) - (s_{1} + l) = s_{2} - s_{1}. 
\end{align*}
By Case~2.2, 
the equality holds in the above, 
and $G_{1} - u_{1}$ is $(s_{1}-1; X)$-degenerate. 
In particular, 
$d_{G_{1}}(u_{1}) = s_{1} + l$.

We now choose such vertices $u_{1}$ and $v_{1}$ so that $d_{G_{1}}(u_{1}) + d_{G_{1}}(v_{1})$ 
is as small as possible. 
Since 
$G_{1} - u_{1}$ is $(s_{1}-1; X)$-degenerate, 
there exist two vertices 
$u_{1}' \in V(G_{1} - u_{1}) \setminus X$ and $v_{1}' \in V(G_{1} - u_{1}) \cap X$ 
with $u_{1}'v_{1}' \notin E(G_{1} - u_{1})$ 
such that $d_{G_{1} - u_{1}}(u_{1}') + d_{G_{1} - u_{1}}(v_{1}') \le 2(s_{1}-1)$ 
(recall that 
by Subclaim~\ref{sub1}, 
$d_{G_{1}-u_{1}}(u) > s_{1} - 1$ for $u \in V(G_{1} - u_{1}) \setminus X$). 
Then by the choice of $u_{1}$ and $v_{1}$, 
\begin{align*}
2s_{1} 
&
= (s_{1} + l) + (s_{1} - l) 
= d_{G_{1}}(u_{1}) + d_{G_{1}}(v_{1}) \\
& 
\le d_{G_{1}}(u_{1}') + d_{G_{1}}(v_{1}') 
= d_{G_{1}-u_{1}}(u_{1}') + d_{G_{1}-u_{1}}(v_{1}') + |E(G) \cap \{u_{1}u_{1}', u_{1}v_{1}'\}| \\
&\le 2(s_{1} - 1) + 2 = 2s_{1}. 
\end{align*}
Thus the equality holds in the above. 
The equality $|E(G) \cap \{u_{1}u_{1}', u_{1}v_{1}'\}| = 2$ implies that 
$v_{1} \neq v_{1}'$ because $u_{1}v_{1} \notin E(G)$. 
Since $d_{G_{1}}(u_{1}) + d_{G_{1}}(v_{1}) = d_{G_{1}}(u_{1}') + d_{G_{1}}(v_{1}')$, 
we can replace $u_{1}$ and $v_{1}$ with $u_{1}'$ and $v_{1}'$, respectively. 
Therefore, Subclaim~\ref{sub2} also holds for the vertex $v_{1}'$, 
and hence 
$v_{1}$ and $v_{1}'$ are two distinct vertices 
such that 
$d_{G}(v_{1}) \le s_{1} + s_{2} - 3$ 
and 
$d_{G}(v_{1}') \le s_{1} + s_{2} - 3$, 
which contradicts (\ref{degree of X}).

This completes the proof of Claim~\ref{Z is non-empty}. 
\qed

\begin{claim}
\label{low degree vertex} 
Suppose that there exists a vertex $u_{2}$ of $G_{2}$ 
such that $zu_{2} \in E(G)$ for all $z \in Z_{1}$. 
Then, the following hold.  
\begin{enumerate}[{\upshape(i)}]
\item There exists a vertex $v_{1}$ of $G_{1}$ 
such that $d_{G}(v_{1}) \le s_{1} + s_{2} - 2$. 
\item If $d_{G_{2}}(u_{2}) \le s_{2} - 1$, 
then 
$d_{G_{1}}(u_{2}) - d_{G_{2}}(u_{2}) \ge s_{1} - s_{2} + 3$. 
\end{enumerate}
\end{claim}
\proof 
(i) Let $z \in Z_{1}$. 
Then by the definition of $Z_{1}$ and Case 2.2(ii), 
$G_{1} - z$ is $(s_{1}-1; X)$-degenerate. 
Suppose first that 
there exists a vertex $u_{1}$ of $V(G_{1} - z) \setminus X$ such that $d_{G_{1} - z}(u_{1}) \le s_{1} - 1$. 
Then, 
$d_{G_{1}}(u_{1}) = d_{G_{1} - z}(u_{1}) + |E(G) \cap \{u_{1}z\}| \le (s_{1} - 1) + 1 = s_{1}$. 
Combining this with Case~2.2(i), 
we also get $d_{G_{2}}(u_{1}) \le s_{2}$, 
and hence $d_{G}(u_{1}) = d_{G_{1}}(u_{1}) + d_{G_{2}}(u_{1}) \le s_{1} + s_{2}$. 
On the other hand, 
it follows from the definition of $X$ that 
$d_{G}(u_{1}) \ge s_{1} + s_{2}$. 
Therefore, the equality holds in the above. 
The equalities $d_{G_{1}}(u_{1}) = s_{1}$ and $d_{G_{2}}(u_{1}) = s_{2}$ imply that $u_{1} \in Z_{1}$. 
The equality $|E(G) \cap \{u_{1}z\}| = 1$ implies that $u_{1}z \in E(G)$. 
But then, from the assumption of Claim~\ref{low degree vertex}, 
$G[\{u_{2}, z, u_{1}\}]$ forms a triangle, a contradiction. 
Thus 
$d_{G_{1} - z}(u) > s_{1} - 1$
for $u \in V(G_{1} - z) \setminus X$.

Since $G_{1}-z $ is $(s_{1}-1; X)$-degenerate, 
this implies that 
there exist two vertices 
$u_{1} \in V(G_{1} - z) \setminus X$ and $v_{1} \in V(G_{1} - z) \cap X$ 
with $u_{1}v_{1} \notin E(G_{1} - z)$ 
such that $d_{G_{1} - z}(u_{1}) + d_{G_{1} - z}(v_{1}) \le 2(s_{1}-1)$. 
Since 
$d_{G_{1} - z}(u_{1}) > s_{1} - 1$, 
we have $d_{G_{1} - z}(v_{1}) \le s_{1} - 2$, 
and hence $d_{G_{1}}(v_{1}) \le s_{1} - 1$. 
Combining this with Case~2.2(i), 
we also get $d_{G_{2}}(v_{1}) \le s_{2} - 1$, 
and hence $d_{G}(v_{1}) \le s_{1} + s_{2} -2$. 
Thus (i) is proved.

(ii) It suffices to show that $d_{G_{1}}(u_{2}) \ge s_{1} + 2$. 
To show it, 
let $z \in Z_{1}$ and $v_{1}$ be a vertex of $G_{1}$ 
such that $d_{G}(v_{1}) \le s_{1} + s_{2} - 2$ 
(we can take such a vertex $v_{1}$ by (i)). 
By the degree condition of $G$, 
$zv_{1} \in E(G)$. 
Recall that by the assumption of Claim~\ref{low degree vertex}, 
$zu_{2} \in E(G)$, 
and hence 
$u_{2}v_{1} \notin E(G)$ because $g(G) \ge 4$. 
Therefore, 
if 
$d_{G_{1}}(u_{2}) \le s_{1} + 1$, 
then 
\begin{align*}
\sigma_{2}(G) 
&\le d_{G}(v_{1}) + d_{G}(u_{2}) 
= d_{G}(v_{1}) + \big( d_{G_{2}}(u_{2}) + d_{G_{1}}(u_{2}) \big)\\
&\le (s_{1} + s_{2} - 2) + \big(  (s_{2} - 1) + (s_{1} + 1) \big) 
= 2(s_{1} + s_{2}) - 2, 
\end{align*}
a contradiction. 
Thus $d_{G_{1}}(u_{2}) \ge s_{1} + 2$. 
\qed

\begin{claim}
\label{good vertices} 
There exist two vertices $z_{1} \in Z_{1}$ 
and $u_{2} \in V(G_{2})$ 
satisfying one of the following {\rm (A)} and {\rm (B)}. 
\begin{enumerate}[{\upshape(A)}]
\item $z_{1}u_{2} \notin E(G)$ and $d_{G_{1}}(u_{2}) - d_{G_{2}}(u_{2}) \ge s_{1} - s_{2} + 1$, 
or 
\item $z_{1}u_{2} \in E(G)$ and $d_{G_{1}}(u_{2}) - d_{G_{2}}(u_{2}) \ge s_{1} - s_{2} + 3$.  
\end{enumerate}
\end{claim}
\proof 
We now apply Lemma~\ref{s_{i} + l + 1} with $\varepsilon = 0$. 
By the symmetry $s_{1}, s_{2}$ in Lemma~\ref{s_{i} + l + 1} 
and the proof, 
we actually have the following 
(see the vertex $u$ of Case~1 
and the vertices $u$ and $v$ of Case~2 in the proof of Lemma~\ref{s_{i} + l + 1}): 
$G_{2}$ satisfies one of (A') and (B'). 
\begin{enumerate}[{\upshape(A')}]
\item 
There exists a vertex $u_{2}$ of $G_{2}$ such that 
$d_{G_{1}}(u_{2}) - d_{G_{2}}(u_{2}) \ge s_{1} - s_{2} + 1$ 
and 
$d_{G_{2}}(u_{2}) \le s_{2} - 1$, or 
\item 
there exist two vertices $u_{2}$ and $v_{2}$ of $G_{2}$ such that 
$d_{G_{1}}(u_{2}) - d_{G_{2}}(u_{2}) \ge s_{1} - s_{2} + 1$ 
and 
$d_{G}(v_{2}) \le s_{1} + s_{2} - 2$. 
\end{enumerate}
We further take $z_{1} \in Z_{1}$, 
and choose it so that $z_{1}u_{2} \notin E(G)$ if possible. 
If $z_{1}u_{2} \notin E(G)$, 
then (A') and (B') directly imply (A). 
Thus we may assume that $z_{1}u_{2} \in E(G)$. 
Then 
by the choice of $z_{1}$, 
we have $zu_{2} \in E(G)$ for all $z \in Z_{1}$. 
It follows from (\ref{degree of X}) and Claim~\ref{low degree vertex}(i) 
that (B') does not occur. 
Thus, (A') holds. 
By Claim~\ref{low degree vertex}(ii), 
the inequality $d_{G_{2}}(u_{2}) \le s_{2} - 1$ implies that 
$d_{G_{1}}(u_{2}) - d_{G_{2}}(u_{2}) \ge s_{1} - s_{2} + 3$. 
Thus (B) holds. 
\qed

Let $z_{1}$ and $u_{2}$ be the same as in Claim~\ref{good vertices}. 
Consider the graphs 
\begin{align*}
\begin{array}{ll}
G_{1}' = G_{1} - z_{1}, &\quad G_{2}' = G[V(G_{2}) \cup \{z_{1}\}],  \\[1mm]
G_{1}'' = G[V(G_{1}') \cup \{u_{2}\}], &\quad G_{2}'' = G_{2}' - u_{2}. 
\end{array}
\end{align*}
Note that $(G_{1}', G_{2}')$ and $(G_{1}'', G_{2}'')$ are partitions of $G$, respectively. 
Then by Lemma~\ref{w(G_{1}, G_{2})}, Claim~\ref{good vertices} and 
the definitions of $Z_{1}, G_{1}', G_{2}', G_{1}''$ and $G_{2}''$, 
we get 
\begin{align}
\label{w(G_{1}'', G_{2}'')}
f(G_{1}'', G_{2}'', s_{1}, s_{2}) \ge f(G_{1}', G_{2}', s_{1}, s_{2}) + 1 
= f(G_{1}, G_{2}, s_{1}, s_{2}) + 1. 
\end{align}

\begin{claim}
\label{property of G_{i}''} 
Each $G_{i}''$ is $(s_{i}; X)$-degenerate. 
\end{claim}
\proof 
By Case~2.2(ii) and the definitions of $Z_{1}$ and $G_{1}'$, 
it follows that 
$G_{1}'$ is $(s_{1}-1; X)$-degenerate, 
and hence by Lemma~\ref{feasible and degenerate}(\ref{(s-1; X)-dege to (s; X)-dege}) 
and the definition of $G_{1}''$, 
we see that 
$G_{1}''$ is $(s_{1}; X)$-degenerate 
(recall that $\omega(G) \le 2 \le s_{1}$). 
Combining this with Lemma~\ref{G-G' is non-complete tri.-free}, 
we also see that 
$G_{2}'' \ (= G - G_{1}'')$ is non-complete. 
Since $G_{2}$ is $(s_{2}-1; X)$-degenerate, 
it follows from Lemma~\ref{feasible and degenerate}(\ref{(s-1; X)-dege to (s; X)-dege}) and the definition of $G_{2}'$ that 
$G_{2}'$ is $(s_{2}; X)$-degenerate. 
Hence, $G_{2}'' \ (= G_{2}' - u_{2})$ is clearly $(s_{2}; X)$-degenerate because $G_{2}''$ is non-complete. 
\qed

By Claim~\ref{property of G_{i}''} and (\ref{w(G_{1}'', G_{2}'')}), 
if 
$G_{1}''$ is $(s_{1}-1; X)$-degenerate, 
or 
if 
$G_{2}''$ is $(s_{2}-1; X)$-degenerate, 
then 
this contradicts (C1). 
Thus, each $G_{i}''$ is not $(s_{i}-1; X)$-degenerate. 
Since each $G_{i}''$ is non-complete, 
this together with Remark~\ref{definition corollary} implies that 
each $G_{i}''$ contains an $(s_{i} - 1 ; X)$-feasible induced subgraph. 
Thus, $(G_{1}'', G_{2}'')$ is the desired partition $(G_{1}^{*}, G_{2}^{*})$.

This completes the proof of Theorem~\ref{feasible pair tri.-free}. 
\qed

\section{Applications}
\label{applications} 

\subsection{Degree conditions for vertex-disjoint cycles} 
\label{deg conditions}

In Sections~\ref{introduction}--\ref{proof of tri.-free}, 
we have considered the existence of disjoint subgraphs with degree conditions, 
and we have shown Theorems~\ref{disjoint subgraphs sigma2} and \ref{disjoint subgraphs sigma2 tri-free}, 
which correspond to Step~1 (the existence of two disjoint subgraphs of high minimum degree sum) 
for Problems~\ref{partition problem sigma2} and \ref{partition problem sigma2 tri-free}. 
These types of results are sometimes useful tools to get degree conditions for packing of graphs, 
i.e., the existence of $k$ disjoint subgraphs 
which belong to some fixed class of graphs. 
In this section, 
we explain it by taking disjoint cycles for example. 
In particular, 
we will give a sharp $\sigma_{2}$ condition for the existence of $k$ disjoint cycles of lengths $0$-mod~$3$ 
by using Theorem~\ref{disjoint subgraphs sigma2 tri-free} 
(see statement (S\ref{s2}), Proposition~\ref{inductive step sigma2} and Theorem~\ref{0modulo3 sigma2 k=1}).

\medskip

In \cite{Chen and Saito}, 
Chen and Saito gave a minimum degree condition for the existence of a cycle of length $0$-mod~$3$. 
Here, 
a cycle $C$ is called a \textit{cycle of length $0$-mod~$3$} 
if $|C| \equiv 0 \pmod 3$.

\begin{Thm}[Chen and Saito \cite{Chen and Saito}] 
\label{Chen and Saito}
Every graph $G$ with $\delta(G) \ge 3$ contains a cycle of length $0$-mod~$3$.
\end{Thm}

As a natural generalization of this result, 
we can consider the following problem.

\begin{problem} 
\label{0modulo3 min deg}
Is the following statement true for any $k \ge 1$? 
\begin{enumerate}
\renewcommand{\labelenumi}{\upshape{(S\arabic{enumi}})}
\setcounter{enumi}{\value{mymemory2}}
\item
\label{s1} 
Every graph $G$ with $\delta(G) \ge 3k$ contains $k$ disjoint cycles of lengths $0$-mod~$3$.
\setcounter{mymemory2}{\value{enumi}}
\end{enumerate}
\end{problem}

In statement (S\ref{s1}), 
the minimum degree condition is best possible if it's true. 
Let $k$ and $n$ be integers with $k \ge 1$  and $n \ge 6k-2$, 
and consider the complete bipartite graph $K_{3k-1, n - 3k + 1}$. 
The minimum degree of this graph is clearly $3k-1$, 
and  
every cycle of length $0$-mod~$3$ in this graph has order at least $6$, 
and hence it does not contain $k$ disjoint cycles of lengths $0$-mod~$3$. 
In addition, 
considering this example, 
we can also consider the more general problem as follows.

\begin{problem} 
\label{0modulo3 sigma2}
Is the following statement true for any $k \ge 1$? 
\begin{enumerate}
\renewcommand{\labelenumi}{\upshape{(S\arabic{enumi}})}
\setcounter{enumi}{\value{mymemory2}}
\item
\label{s2} 
Every graph $G$ of order at least $3k$ with 
$\sigma_{2}(G) \ge 6k - 1$ contains $k$ disjoint cycles of lengths $0$-mod~$3$. 
\setcounter{mymemory2}{\value{enumi}}
\end{enumerate}
\end{problem}

Since $\sigma_{2}(K_{3k-1, n - 3k + 1}) = 6k-2$, 
the graph 
$K_{3k-1, n - 3k + 1}$ shows that 
``$\sigma_{2}(G) \ge 6k-1$'' cannot be replaced by ``$\sigma_{2}(G) \ge 6k-2$'' in statement (S\ref{s2}). 
Moreover, 
since $\sigma_{2}(G) \ge 2\delta(G)$ for a graph $G$, 
it follows that 
statement (S\ref{s2}) is stronger than statement (S\ref{s1}) 
(note that $|G| \ge 3k$ in statement (S\ref{s1})).

In order to attack the above problems, 
one may use induction on $k$. 
In particular, 
for Problem~\ref{0modulo3 min deg}, 
we already know that 
statement (S\ref{s1}) is true when $k=1$ by Theorem~\ref{Chen and Saito},  
that is, Problem~\ref{0modulo3 min deg} can be solved by showing the inductive step. 
In the inductive step of this type of problem, 
Theorems~\ref{Stiebitz}, \ref{Kaneko}, \ref{disjoint subgraphs sigma2} and \ref{disjoint subgraphs sigma2 tri-free} sometimes can work effectively. 
In fact, 
we can easily obtain the following 
by using Theorems~\ref{Kaneko} and \ref{disjoint subgraphs sigma2 tri-free}, respectively.

\begin{prop} 
\label{inductive step min deg}
If statement {\rm (S\ref{s1})} is true for $k = 1$, 
then statement {\rm (S\ref{s1})} is true for any $k \ge 1$. 
\end{prop}

\begin{prop} 
\label{inductive step sigma2}
If statement {\rm (S\ref{s2})} is true for $k=1$, 
then statement {\rm (S\ref{s2})} is true for any $k \ge 1$. 
\end{prop}

We only prove Proposition~\ref{inductive step sigma2} 
because we can obtain Proposition~\ref{inductive step min deg} 
by the same argument.

\medskip
\noindent
\textbf{Proof of Proposition~\ref{inductive step sigma2}.}~We show that 
statement (S\ref{s2}) is true for any $k \ge 1$ 
by induction on $k$. 
By the assumption of Proposition~\ref{inductive step sigma2}, 
statement (S\ref{s2}) is true when $k=1$. 
Suppose that statement (S\ref{s2}) is true up to the row $k-1$, $k\geq 2$, and let us study for $k$.

Let $G$ be a graph of order at least $3k$ with 
$\sigma_{2}(G) \ge 6k - 1$. 
We show that 
$G$ contains $k$ disjoint cycles of lengths $0$-mod~$3$. 
If $G$ is complete, 
then the assertion clearly holds. 
Thus we may assume that $G$ is non-complete.

Suppose first that $G$ contains a triangle $C$. 
Then every vertex of $G$ not in $C$ has at most $3$ neighbors in $C$, 
and hence 
$\sigma_{2}(G - C) \ge (6k-1) - 6 = 6(k-1) - 1$. 
Note that $|G-C| \ge 3(k-1)$. 
Since 
statement (S\ref{s2}) is true for $k-1$ by the induction hypothesis, 
$G-C$ contains $k-1$ disjoint cycles of lengths $0$-mod~$3$. 
With the cycle $C$, we get then $k$ disjoint cycles of lengths $0$-mod~$3$ in $G$.

Suppose now that $g(G) \ge 4$. 
Then, since 
$\sigma_{2}(G) \ge 6k - 1 = 2\big( 3(k-1) + 3 \big) - 1$, 
it follows from Theorem~\ref{disjoint subgraphs sigma2 tri-free} 
that 
there exist two disjoint subgraphs $H_{1}$ and $H_{2}$ of $G$ 
such that 
$|H_{1}| \ge 2 \cdot 3(k-1) > 3(k-1)$, 
$\sigma_{2}(H_{1}) \ge 2 \cdot 3(k-1) - 1 = 6(k-1) - 1$, 
$|H_{2}| \ge 2 \cdot 3 > 3$ and 
$\sigma_{2}(H_{2}) \ge 2 \cdot 3 - 1 = 5$. 
Hence by the induction hypothesis, 
$H_{1}$ contains $k-1$ disjoint cycles of lengths $0$-mod~$3$, 
and $H_{2}$ contains a cycle of length $0$-mod~$3$. 
We get then $k$ disjoint cycles of lengths $0$-mod~$3$ in $G$. 
\qed

By Theorem~\ref{Chen and Saito} and Proposition~\ref{inductive step min deg}, 
we see that Problem~\ref{0modulo3 min deg} is solved in affirmative. 
Similarly, 
by Proposition~\ref{inductive step sigma2}, 
it is only necessary to consider the case of $k=1$ for Problem~\ref{0modulo3 sigma2}. 
In the next subsection, 
we completely solve Problem~\ref{0modulo3 sigma2} by showing the following.

\begin{thm} 
\label{0modulo3 sigma2 k=1}
Every graph $G$ of order at least $3$ with $\sigma_{2}(G) \ge 5$ contains a cycle 
of length $0$-mod~$3$.
\end{thm}

\medskip

We mention other cases in the rest of this subsection. 
It is well known that 
every graph $G$ of order at least $3k$ with $\delta(G) \ge 2k$ contains 
$k$ disjoint cycles, 
which is a classical result by Corr\'{a}di and Hajnal \cite{Corradi-Hajnal}. 
By a similar argument as in the proof of Proposition~\ref{inductive step sigma2}, 
we can easily obtain a slightly weaker version of this: 
Every graph $G$ with $\delta(G) \ge 3k-1$ contains 
$k$ disjoint cycles. 
Because if $G$ is a graph with $\delta(G) \ge 3k-1$, 
then Theorem~\ref{Stiebitz} implies that 
$G$ contains disjoint subgraphs $H_{1}$ and $H_{2}$ 
such that 
$\delta(H_{1}) \ge 3(k-1) - 1$ and $\delta(H_{2}) \ge 2$ 
(see also \cite[Corollary 2]{Stiebitz}). 
Similarly, 
Theorem~\ref{Kaneko} leads to a triangle-free version of 
Corr\'{a}di and Hajnal's theorem 
(note that the minimum degree condition is best possible even if we assume a triangle-freeness). 
In addition, 
Enomoto \cite{Enomoto} and Wang \cite{Wang} independently gave a $\sigma_{2}$-version as follows: 
Every graph $G$ of order at least $3k$ with $\sigma_{2}(G) \ge 4k-1$ contains 
$k$ disjoint cycles. 
Theorems~\ref{disjoint subgraphs sigma2} and \ref{disjoint subgraphs sigma2 tri-free} 
immediately lead to slightly weaker versions of this result as above.

Recently, 
Gould, Horn and Magnant \cite{Gould-Horn-Magnant} 
proposed the following conjecture, 
which is a common generalization 
of the above Corr\'{a}di and Hajnal's theorem and Hajnal and Szemer\'{e}di's theorem \cite{Hajnal and Szemeredi} 
``every graph $G$ of order exactly $(c+1)k$ and of minimum degree at least $ck$ contains $k$ disjoint complete graphs of orders $c+1$''. 
Here, for an integer $c \ge 0$, 
a $c$-\textit{chorded cycle} is a cycle with $c$ chords.

\begin{Conj}[Gould, Horn and Magnant \cite{Gould-Horn-Magnant}]
\label{Gould-Horn-Magnant conj}
Let $c, k$ be integers with $c \ge 2$ and $k \ge 1$. 
Every graph $G$ of order at least $(c+1)k$ with $\delta(G) \ge ck$ 
contains $k$ disjoint $\frac{(c+1)(c-2)}{2}$-chorded cycles. 
\end{Conj}

They showed that 
this conjecture is true 
for very large graphs compared to $c$ and $k$ 
(see \cite{Gould-Horn-Magnant} for more details). 
However, 
by using Theorem~\ref{Stiebitz} and the following result (Theorem~\ref{Gould-Horn-Magnant}) by Gould et al. \cite{Gould-Horn-Magnant}, 
we can easily obtain a slightly weaker version of Conjecture~\ref{Gould-Horn-Magnant conj} 
(see Corollary~\ref{weak Gould-Horn-Magnant conj}).

\begin{Thm}[Gould, Horn and Magnant \cite{Gould-Horn-Magnant}]
\label{Gould-Horn-Magnant}
Let $c \ge 2$ be an integer. 
Every graph $G$ with $\delta(G) \ge c$ contains a $\frac{(c+1)(c-2)}{2}$-chorded cycle. 
\end{Thm}

\begin{cor}
\label{weak Gould-Horn-Magnant conj}
Let $c, k$ be integers with $c \ge 2$ and $k \ge 1$. 
Every graph $G$ with $\delta(G) \ge (c+1)k - 1$ 
contains $k$ disjoint $\frac{(c+1)(c-2)}{2}$-chorded cycles. 
\end{cor}

In the same paper, 
they also showed that 
Conjecture~\ref{Gould-Horn-Magnant conj} holds when $G$ is triangle-free (see \cite[Theorem 6]{Gould-Horn-Magnant}).
However, 
we can also obtain it by using Theorems~\ref{Kaneko} and \ref{Gould-Horn-Magnant}. 
Moreover, 
if we can obtain a $\sigma_{2}$-version of Theorem~\ref{Gould-Horn-Magnant}, 
then 
by combining it with 
Theorem~\ref{disjoint subgraphs sigma2}, 
we can also get the $\sigma_{2}$-version of Corollary~\ref{weak Gould-Horn-Magnant conj}.

\subsection{Proof of Theorem~\ref{0modulo3 sigma2 k=1}} 
\label{proof of 0modulo3 sigma2 k=1}

In this subsection, 
we prove Theorem~\ref{0modulo3 sigma2 k=1}. 
In the proof, 
we will use the following 
result of Chen and Saito which is stronger than Theorem~\ref{Chen and Saito}.

\begin{Thm}[Chen and Saito \cite{Chen and Saito}] 
\label{Chen and Saito2}
Every graph $G$ of order at least $3$ with at most one vertex of degree less than $3$, contains a cycle of length $0$-mod~$3$. 
\end{Thm}

For a graph $G$ of order $n$ and an unordered pair $\lbrace u,v\rbrace$ of distinct vertices of $G$ (adjacent or not), 
we define the graph $G_{u,v}$ of order $n-1$, as follows:
\begin{itemize}
\item The vertices of $G_{u,v}$ are the vertices $x$ of $G$ distinct from $u$ and $v$, and the pair $\lbrace u,v\rbrace$. 
\item The edges of $G_{u,v}$ are the edges $xy$ of $G$ with $x, y \notin \lbrace u,v\rbrace$ 
and the edges $xy$ with $x = \lbrace u,v\rbrace$, $y \notin \lbrace u,v\rbrace$ 
and $E(G) \cap \{yu, yv\} \neq \emptyset$. 
\end{itemize}

\medskip

\noindent
For a graph $G$ and an integer $s$, 
we further let 
$V_{s}(G) = \{v \in V(G) : d_{G}(v) = s\}$ 
and 
$V_{ \ge s}(G) = \{v \in V(G) : d_{G}(v) \ge s\}$.

\medskip

Now we are ready to prove Theorem~\ref{0modulo3 sigma2 k=1}.

\medskip
\noindent
\textit{Proof of Theorem~\ref{0modulo3 sigma2 k=1}.}~We proceed by induction on $n := |G|$. 
Clearly the assertion is true for $n = 3$. 
Suppose that the assertion is true up to the row $n-1$, $n\geq 4$, and let us study for $n$. 
So $G$ is a graph of order $n\geq 4$ with $\sigma_{2}(G) \ge 5$.
Clearly, if $G$ contains triangles, we are done. 
So, we may suppose that $g(G) \ge 4$. 
Since $\sigma_{2}(G) \ge 5$ and $g(G) \ge 4$, 
it follows from Lemma~\ref{sigma}(\ref{complete}) that 
$V_{\le 2}(G)$ is a clique and $|V_{\le 2}(G)| \le 2$. 
If $|V_{\le 2}(G)| \le 1$, by Theorem~\ref{Chen and Saito2}, $G$ contains a cycle of length $0$-mod~$3$.
So, we may suppose that $|V_{\le 2}(G)| = 2$, say $V_{\le 2}(G) = \lbrace x, y \rbrace$. 
Recall that $xy \in E(G)$. 
Suppose that one of the vertices of $V_{\le 2}(G)$, say $y$ is of degree $1$ in $G$. 
It is easy to see that the induced subgraph $G_1=G-y$ of $G$ 
has exactly one vertex of degree less than $3$. 
Then $G_1$ contains a cycle of length $0$-mod~$3$, and we are done. 
So $V_{\le 2}(G) = \{x, y\} = V_{2}(G)$,  
and since $x$ and $y$ does not have a common neighbor (for otherwise triangles), 
$x$ has a unique neighbor $u$ distinct from $y$ and $y$ has a unique neighbor $v$ distinct from $x$ and $u$. 
Suppose first that $uv \in E(G)$. We put $G'=G-\lbrace x,y\rbrace$. 
It is easy to see that every vertex of $G'_{u,v}$ distinct from the vertex $\lbrace u,v\rbrace$ is of degree at least $3$ in $G'_{u,v}$. 
Then $G'_{u,v}$ contains a cycle $C$ of length $0$-mod~$3$. 
If $\lbrace u,v\rbrace$ is not a vertex of $C$, clearly we are done. 
Suppose now that $\{u, v\}$ is a vertex of $C$. 
We put $C=(x_1, x_2,\ldots,x_r,x_1)$, where $x_{1} = \{u, v\}$. 
If $u$ or $v$ is a common neighbor of $x_{r}$ and $x_{2}$, clearly we are done. 
If it is not the case, we may suppose that $x_2$ is adjacent to $v$ and that $x_r$ is adjacent to $u$. 
Then $P=(v,x_2,\ldots,x_r,u)$ is a path of $G$ such that $|P| \equiv 1 \pmod 3$. 
Since $x, y \notin V(P)$, 
$C'=(x,y,P,x)$ is a cycle of $G$ of length $0$-mod~$3$.

\vspace{0.2cm} 
So, we may suppose that $uv \notin E(G)$. 
We distinguish two cases.

\begin{enumerate}[{\textup{{\bf Case~\arabic{enumi}.}}}]
\setcounter{enumi}{0}
\item $\{u, v\} \cap V_{\ge 4}(G) \neq \emptyset$. 
\end{enumerate}

Without loss of generality, we may suppose that $u \in V_{\ge 4}(G)$. 
Then it is easy to see that every vertex of $G'=G-\lbrace x,y\rbrace$ distinct from $v$ is of degree at least $3$ in $G'$. 
By Theorem~\ref{Chen and Saito2}, 
$G'$ (and therefore $G$) contains a cycle of length $0$-mod~$3$, and then we are done.

\begin{enumerate}[{\textup{{\bf Case~\arabic{enumi}.}}}]
\setcounter{enumi}{1}
\item $\{u, v\} \cap V_{\ge 4}(G) = \emptyset$, i.e., $u, v \in V_{3}(G)$. 
\end{enumerate}

We put again $G'=G-\lbrace x,y\rbrace$, and we consider the graph $G'_{u,v}$. 
Suppose that the vertex $\lbrace u,v\rbrace$ of $G'_{u,v}$ has degree at least $3$ in $G'_{u,v}$. 
Since each of the vertices $u$ and $v$ has exactly $2$ neighbors in $G-\lbrace x,y,v,u\rbrace$, 
and since $\lbrace u,v\rbrace$ has degree at least $3$ in $G'_{u,v}$, 
it follows that $u$ and $v$ have at most one common neighbor in $G'$, 
and then it is easy to see that at most one vertex of $G'_{u,v}$ is of degree less than $3$ in $G'_{u,v}$. 
By Theorem~\ref{Chen and Saito2}, $G'_{u,v}$ contains a cycle $C$ of length $0$-mod~$3$, 
and then as above we get a cycle of $G$ of length $0$-mod~$3$. 
So, we may suppose that $\lbrace u,v\rbrace$ has degree at most $2$ in $G'_{u,v}$. 
Then necessarily, $u$ and $v$ have two common neighbors $w$ and $z$ in $G'$. 
Observe that $x, y, v $ and $u$ have no neighbors in $G-\lbrace x,y,v,u, w,z\rbrace$ and that $wz \notin E(G)$.

Suppose that $\{w, z\} \cap V_{\ge 4}(G) \neq \emptyset$. 
Consider then the graph $G''=G-\lbrace x,y,v\rbrace$. 
Then, it is easy see that $\sigma_2(G'')\geq 5$, and hence by the induction hypothesis, we are done. 
So, we may suppose $w, z \in V_{3}(G)$.  
We consider now two subcases.

\begin{enumerate}[{\textup{{\bf Case~2.\arabic{enumi}.}}}]
\setcounter{enumi}{0}
\item $w$ has a neighbor $a$ in $G-\lbrace x,y,v,u, w,z\rbrace$ 
and $z$ has a neighbor $b$ in $G-\lbrace x,y,v,u, w,z\rbrace$ distinct from $a$. 
\end{enumerate}

Suppose first that $ab \in E(G)$. 
We consider the graph $G_1=G-\lbrace x,y,v,u, w,z \rbrace$. 
It is easy to see that $\sigma_2(G_1)\geq 5$ and then by the induction hypothesis we are done. 
Suppose now that $ab \notin E(G)$. 
Then the graph $G_2=G_1+ab$ is of minimum degree at least $3$. 
Then $G_2$ contains a cycle $C_1$ of length $0$-mod~$3$. 
If $C_1$ does not contain the edge $ab$ of $G_2$, we are done. 
If $C_1$ contains $ab$, then by deleting this edge and by adding the vertices $w$, $u$ and $z$, 
we get a cycle of $G$ of length $0$-mod~$3$, and so we are done.

\begin{enumerate}[{\textup{{\bf Case~2.\arabic{enumi}.}}}]
\setcounter{enumi}{1}
\item The vertices $w$ and $z$ have a common neighbor $a$ in $G-\lbrace x,y,v,u,w,z\rbrace$. 
\end{enumerate}

It is easy to see that all the vertices of $G_1=G-\lbrace x,y,v,u, w,z\rbrace$ distinct from $a$ are of degree at least $3$ in $G_1$, 
and then by Theorem~\ref{Chen and Saito2}, we are done. 
So, the assertion is true for $n$, and Theorem~\ref{0modulo3 sigma2 k=1} is proved. 
\qed





\begin{thebibliography}{}

\bibitem {Bazgan et al.}
C.~Bazgan, Z.~Tuza, D.~Vanderpooten, 
\textit{Efficient algorithms for decomposing graphs under degree constraints}, 
Discrete Appl. Math. \textbf{155} (2007) 979--988.  


\bibitem {Chen and Saito}
G.~Chen, A.~Saito,  
\textit{Graphs with a cycle of length divisible by three}, 
J. Combin. Theory Ser. B \textbf{60} (1994) 277--292.  


\bibitem{Corradi-Hajnal} 
K.~Corr\'adi, A.~Hajnal, 
\textit{On the maximal number of independent circuits in a graph}, 
Acta Math. Acad. Sci. Hungar. \textbf{14} (1963) 423--443.


\bibitem {Diestel}
R.~Diestel, 
\textit{Graph Theory}, Fourth edition. Graduate Texts in Mathematics, 173, Springer, Heidelberg, 2010.


\bibitem {Enomoto} 
H.~Enomoto, 
\textit{On the existence of disjoint cycles in a graph}, 
Combinatorica \textbf{18} (1998) 487--492.  


\bibitem{Diwan} 
A.~Diwan, 
\textit{Decomposing graphs with girth at least five under degree constraints}, 
J. Graph Theory \textbf{33} (2000) 237--239.


\bibitem{Dunbar-Frick}
J.E.~Dunbar, M.~Frick, 
\textit{The Path Partition Conjecture is true for claw-free graphs}, 
Discrete Math. \textbf{307} (2007) 1285--1290. 


\bibitem{Dunbar-Frick-Bullock}
J.E.~Dunbar, M.~Frick, F.~Bullock, 
\textit{Path partitions and $P_{n}$-free sets}, 
Discrete Math. \textbf{289} (2004) 145--155.



\bibitem{Gould-Horn-Magnant}
R.~Gould, P.~Horn, C.~Magnant, 
\textit{Multiply Chorded Cycles}, 
SIAM Journal on Discrete Math. \textbf{28} (2014) 160--172. 


\bibitem{Hajnal and Szemeredi}
A.~Hajnal, E.~Szemer\'{e}di, 
\textit{Proof of a conjecture of P. Erd\H{o}s}, 
Combinatorial Theory and Its Application \textbf{2} (1970) 601--623. 


\bibitem{Kaneko} 
A.~Kaneko, 
\textit{On decomposition of triangle-free graphs under degree constraints}, 
J. Graph Theory \textbf{27} (1998) 7--9.


\bibitem{Kuhn-Osthus}
D.~K\"{u}hn, D.~Osthus, 
\textit{Partitions of graphs with high minimum degree or connectivity}, 
J. Combin. Theory Ser. B \textbf{88} (2003) 29--43.  


\bibitem{Lovasz} 
L.~Lov\'{a}sz, 
\textit{On decomposition of graphs}, 
Studia Sci. Math. Hungar. \textbf{1} (1966) 237--238.


\bibitem{Stiebitz} 
M.~Stiebitz, 
\textit{Decomposing graphs under degree constraints}, 
J. Graph Theory \textbf{23} (1996) 321--324.


\bibitem{Thomassen1}
C.~Thomassen, 
\textit{Non-separating cycles in $k$-connected graphs}, 
J. Graph Theory \textbf{5} (1981) 351--354.


\bibitem{Thomassen2}
C.~Thomassen, 
\textit{Paths, circuits and subdivisions}, 
in (L. W. Beineke and R.J. Wilson, Eds.), 
Selected topics in graph theory III, academic Press, New York (1988) 97--131.


\bibitem{Thomassen3}
C.~Thomassen, 
\textit{Configurations in graphs of large minimum degree, connectivity, or chromatic number}, 
Combinatorial Mathematics: Proceedings of the Third International Conference (New York, 1985), 402--412,
Ann. New York Acad. Sci., 555, New York Acad. Sci., New York, 1989.


\bibitem {Wang} 
H.~Wang, 
\textit{On the maximum number of independent cycles in a graph}, 
Discrete Math. \textbf{205} (1999) 183--190.


\bibitem{Yang-Vumar}
F.~Yang, E.~Vumar, 
\textit{A note on a cycle partition problem}, 
Appl. Math. Lett. \textbf{24} (2011) 1181--1184. 


\end{thebibliography}
\end{document}